\documentclass[12pt]{amsart}

\newtheorem{theorem}{Theorem}[section]
\newtheorem{lemma}[theorem]{Lemma}
\newtheorem{corollary}[theorem]{Corollary}
\newtheorem{proposition}[theorem]{Proposition}
\newtheorem{remark}[theorem]{Remark}
\newtheorem{exercise}[theorem]{Exercise}
\newtheorem{definition}[theorem]{Definition}

\newcommand{\ncom}{\newcommand}
\ncom{\ep}{\epsilon}
\ncom{\rar}{\rightarrow}
\ncom{\lrar}{\longrightarrow}
\ncom{\ov}{\overline}
\ncom{\m}{\mbox}
\ncom{\sta}{\stackrel}
\ncom{\comx}{{\mathbb C}}
\ncom{\A}{{\mathbb A}}
\ncom{\Z}{{\mathbb Z}}
\ncom{\Q}{{\mathbb Q}}
\ncom{\R}{{\mathbb R}}
\ncom{\G}{{\mathbb G}}
\ncom{\hH}{{\mathbb H}}
\ncom{\al}{\alpha}
\ncom{\p}{{\mathbb P}}
\ncom{\E}{{\mathbb E}}
\ncom{\N}{{\mathbb N}}
\ncom{\K}{{\mathbb K}}
\ncom{\X}{{\mathbb X}}
\ncom{\f}{\frac}
\ncom{\cA}{{\mathcal A}}
\ncom{\cB}{{\mathcal B}}
\ncom{\cD}{{\mathcal D}}
\ncom{\cX}{{\mathcal X}}
\ncom{\cO}{{\mathcal O}}
\ncom{\cW}{{\mathcal W}}
\ncom{\cL}{{\mathcal L}}
\ncom{\cP}{{\mathcal P}}
\ncom{\cH}{{\mathcal H}}
\ncom{\cS}{{\mathcal S}}
\ncom{\cM}{{\mathcal M}}
\ncom{\cC}{{\mathcal C}}
\ncom{\cT}{{\mathcal T}}
\ncom{\cF}{{\mathcal F}}
\ncom{\cN}{{\mathcal N}}
\ncom{\cJ}{{\mathcal J}}
\ncom{\cV}{{\mathcal V}}
\ncom{\cZ}{{\mathcal Z}}
\ncom{\cU}{{\mathcal U}}
\ncom{\cSU}{{\mathcal S \mathcal U}}
\ncom{\cG}{{\mathcal G}}
\ncom{\cQ}{{\mathcal Q}}
\ncom{\cR}{{\mathcal R}}
\ncom{\cY}{{\mathcal Y}}
\ncom{\cE}{{\mathcal E}}

\begin{document}
\baselineskip=16pt


\title[Parabolic Reznikov theorem]{The Chern character of a parabolic bundle, and a parabolic Reznikov theorem 
in the case of finite order at infinity}


\author[J. N. Iyer]{Jaya NN  Iyer}

\address{School of Mathematics, Institute for Advanced Study, 1 Einstein Drive, Princeton NJ 08540 USA.}
\email{jniyer@ias.edu}

\address{The Institute of Mathematical Sciences, CIT
Campus, Taramani, Chennai 600113, India}
\email{jniyer@imsc.res.in}

\author[C. T.  Simpson]{Carlos T  Simpson}
\address{CNRS, Laboratoire J.-A.Dieudonn\'e, Universit\'e de Nice--Sophia Antipolis,
Parc Valrose, 06108 Nice Cedex 02, France}
\email{carlos@math.unice.fr}

\footnotetext{Mathematics Classification Number: 14C25, 14D05, 14D20, 14D21 }

\begin{abstract}
In this paper, we obtain an explicit formula for the Chern character of a locally abelian parabolic bundle in terms of its constituent bundles.
Several features and variants of parabolic structures are discussed. Parabolic bundles arising from logarithmic connections form an important class of examples. As
an application, we consider the situation when the local monodromies are semi-simple and are of finite order at infinity. In this case the parabolic
Chern classes of the associated locally abelian parabolic bundle are deduced to be zero in the rational Deligne cohomology in degrees $\geq 2$. 
\end{abstract}

\maketitle


\section{Introduction}


Parabolic bundles were introduced by Mehta and Seshadri \cite{MehtaSeshadri} \cite{Seshadri} over curves and the definition was extended over higher dimensional 
varieties by Maruyama and Yokogawa \cite{MaruyamaYokogawa} Biswas \cite{Biswas}, Li \cite{Li}, Steer-Wren \cite{Steer-Wren}, Panov \cite{Panov} and Mochizuki \cite{Mochizuki}.
A \textit{parabolic bundle} $F$ on a variety $X$ is a collection of vector bundles $F_\alpha$, indexed by a set of 
\textit{weights}, i.e., $\alpha$ runs over a multi-indexing set $\frac{1}{n}\Z \times \frac{1}{n}\Z \times...\times \frac{1}{n}\Z$, 
for some denominator $n$. Further, all the 
bundles $F_\alpha$ restrict on the complement $X-D$ of some normal crossing divisor $D=D_1+...+D_m$ to the same 
bundle, the index $\alpha$ is an $m$-tuple and the $F_\alpha$ satisfy certain normalization/support hypothesis (see \S{2.1}).

This work is a sequel to \cite{Iy-Si}, which in turn was motivated by Reznikov's work on characteristic classes of
flat bundles \cite{Reznikov}, \cite{Reznikov2}. As a long-range goal we would like to approach the Esnault conjecture \cite{Esnault2} that the
Chern classes of Deligne canonical extensions of motivic flat bundles vanish in the rational Chow groups. Reznikov's
work shows the vanishing of an important piece of these classes, over the subset of definition of a flat bundle.
We think that it should be possible to define secondary classes over a completed variety
for flat connections which are quasi-unipotent
at infinity, and to extend Reznikov's results to this case. At the end of this paper we treat a first and essentially
easy case, when the monodromy transformations at infinity have finite order.
We hope to treat the general case 
in the future and regain an understanding of characteristic classes such as Sasha Reznikov had. 

A different method for obtaining a very partial 
result on the Esnault conjecture, removing a hypothesis from the GRR formula of Esnault-Viehweg
\cite{EsnaultViehweg}, was done in \cite{Iy-Si}. 
There we used a definition of the Chern character obtained
from the correspondence between locally abelian parabolic bundles and usual vector bundles on
a particular Deligne--Mumford stack denoted by $Z_m=  X\langle \frac{D_1}{n},\ldots ,\frac{D_m}{n}\rangle $ (see \cite{Borne}, \cite[\S{2.3}]{Iy-Si}, 
\cite{Cadman}, \cite{MatsukiOlsson}, \cite{MehtaSeshadri}, \cite{Boden}, \cite{Biswas}).
The Chern character of $F$ is defined to be the Chern character of the corresponding vector bundle on this stack. 
This was sufficient for our application in \cite{Iy-Si}, however it is clearly unsatisfactory to have only an
abstract definition rather than a formula. 

The aim of this note is to give an explicit formula for the Chern character in terms of the Chern character of the constituent bundles 
$F_\alpha$ and the divisor components $D_i$ in the rational Chow groups of $X$. This procedure, using a DM stack to define the Chern character
and then giving a computation, was first done for the parabolic degree by Borne in \cite{Borne}, however his techniques are different from ours. 
The parabolic aspect of the problem of extending characteristic classes for bundles from an open variety to its completion
should in the future form a small
part of a generalization of Reznikov's work and we hope the present paper can contribute in that direction. 

With our fixed denominator $n$, introduce the notation
$$
[a_1,\ldots , a_m] := (\frac{a_1}{n}, \ldots , \frac{a_m}{n})
$$
for multi-indices, so the parabolic structure is determined by the bundles
$F_{[a_1,\ldots , a_m]}$ for $0\leq a_i < n$ with $a_i$ integers. 

We prove the following statement.

\begin{theorem}
\label{mainformula}
Suppose $F$ is a locally abelian parabolic bundle on $X$ with respect to $D_1,...,D_m$, with $n$ as the denominator.
Then we have the following
formula for the Chern character of $F$:
\begin{equation}\label{introformula}
{\rm ch}(F) = \frac{\sum _{a_1=0}^{n-1}\cdots \sum _{a_m=0}^{n-1} e^{-\sum _{i=1}^m \f{a_i}{n}D_i}
{\rm ch}(F_{[a_1,\ldots , a_m]})}
{\sum _{a_1=0}^{n-1}\cdots \sum _{a_m=0}^{n-1} e^{-\sum _{i=1}^m \f{a_i}{n}D_i} }.
\end{equation}
In other words, the Chern character of $F$ is the weighted average of the Chern characters of 
the component bundles, with weights $e^{-\sum _{i=1}^m \f{a_i}{n}D_i}$.
\end{theorem}

The proof is by showing that the parabolic bundle obtained by twisting $F$ by a direct sum of line bundles involving 
$D_i$ is \textit{componentwise isomorphic} to a direct sum of the constituent bundles $F_{[a_1,\ldots,a_m]}$ twisted by parabolic line
bundles involving $D_i$ (see Corollary \ref{cchiso}).   
The proof is concluded by proving the main theorem on the invariance of the Chern character 
under componentwise Chow isomorphism (see Theorem \ref{invariance-th}). It says: if locally abelian parabolic bundles 
$F$ and $G$ whose constituent bundles $F_{[a_1,\ldots , a_m]}$  and $G_{[a_1,\ldots , a_m]}$ 
have the same Chern character, for all $a_i$ with $0\leq a_i <n$, then $F$ and $G$ also have the same Chern character in the rational Chow groups
of $X$.

We also give variants of the Chern character formula. One can associate a parabolic structure $F$ to a 
vector bundle $E$ on $X$ and given filtrations on the restriction of $E$ on the divisor components of $D$ (see \S 2, \S 6). 
If $X$ is a surface this is automatically locally abelian, but in higher dimensions it is not 
always the case (see Lemma \ref{surface}).
When the structure is locally abelian, we obtain a formula for ${\rm ch}(F)$ which involves ${\rm ch}(E)$ and terms under the Gysin maps on the multiple 
intersections of the divisor components of $D$ (see Corollary \ref{co.-gysinformula} and Corollary \ref{co.-gysinformulagraded}). 
The shape of the formula depends on the way the 
filtrations intersect on the multiple intersections of the divisor components.

In \S 6 we give two easy counterexamples which show that the Chern character of a parabolic bundle cannot be obtained easily from
just the Chern character of the underlying bundle and that of its filtrations taken separately, nor from the data of a filtration of
subsheaves indexed by a single parameter for the whole divisor (Maruyama-Yokogawa's original definition \cite{MaruyamaYokogawa}). These show that in order to
obtain a good formula we should consider all of the bundles $F_{[a_1,\ldots , a_m]}$. This version of parabolic structure was
first introduced by Li \cite{Li}, Steer-Wren \cite{Steer-Wren} and Mochizuki \cite{Mochizuki}. 

We treat parabolic bundles with real weights in \S 8. The aim is to define pullback of a locally 
abelian parabolic bundle as a locally abelian parabolic bundle. This is done by approximating with 
the rational weights case (see Lemma 8.5). Properties like functoriality, additivity and 
multiplicativity of the Chern character are also discussed. In \S 9, on a smooth surface, parabolic structures at 
multiple points are discussed and a Chern character formula is obtained. Logarithmic connections were 
discussed by Deligne in \cite{Deligne}. We discuss some filtrations defined by the residue transformations 
of the connection at infinity. When the eigenvalues of the residues are rational and non-zero, a 
locally abelian parabolic bundle was associated in \cite{Iy-Si}, and this construction is considered further in \S \ref{reznikov}. 
When the residues are nilpotent, we continue in \S 9 with something different: assign arbitrary weights 
to the pieces of the monodromy weight filtration of the nilpotent residue operators, creating a family of parabolic bundles
indexed by the choices of weights. If $X$ is a surface then these are automatically locally abelian, 
and as an example we make explicit the computation of the parabolic 
Chern character $\m{ch}(F)$ in the case of a weight one unipotent Gauss-Manin system $F$,
see Lemma \ref{weightone}.  

In \S \ref{reznikov} we consider the extension of Reznikov's theory to flat bundles with finite order monodromy at infinity.
Such bundles may be considered as flat bundles over a DM-stack of the form $Z_m=  X\langle \frac{D_1}{n},\ldots ,\frac{D_m}{n}\rangle $,
and Reznikov's theorem \cite{Reznikov2} applies directly (or alternatively, over a finite Kawamata covering). 
The only knowledge which we can add is that our formula of
Theorem \ref{mainformula} gives parabolic Chern classes in terms of the parabolic structure on $X$ deduced from the flat bundle, 
and Reznikov's theorem
can be stated as vanishing of these classes. This might have computational content in explicit examples. 

\begin{proposition}$(\rm{Parabolic\,\, Reznikov\,\, theorem})$:
Suppose $(E_U,\nabla _U)$ is a flat bundle on $U$ with rational and semisimple residues, or equivalently
the monodromy transformations at infinity are of finite order. Let $F$ denote the corresponding locally 
abelian parabolic bundle. Define the Deligne Chern character ${\rm ch}^{\cD}(F)$ by the formula
\eqref{introformula} in the rational Deligne cohomology, and the Chern classes ${\rm c}^{\cD}_p(F)$ by the usual formula. 
Then the classes ${\rm c}^{\cD}_p(F)$ for all $p\geq 2$ vanish.  
\end{proposition}

{\Small
\textbf{Acknowledgements}: We thank P. Deligne for having useful discussions. The first named author is supported by NSF. This material is based upon work supported by the 
National Science Foundation under agreement No. DMS-0111298. Any opinions, findings and conclusions or recommendations 
expressed in this material are those of the authors and do not necessarliy reflect the views of the National Science Foundation. 
}

\section{Parabolic bundles}


Let $X$ be a smooth projective variety over an algebraically closed field of characteristic zero, 
with $D$ a normal crossing divisor on $X$. Write
$D=D_1+\ldots +D_m$ where $D_i$ are the irreducible smooth components and meeting transversally.
We use an approach to parabolic bundles based on multi-indices $(\alpha _1,\ldots , \alpha _m)$
of length equal to the number of components of the divisor. This approach, having its origins in the original
paper of Mehta and Seshadri \cite{MehtaSeshadri},  was introduced in higher dimensions by Li \cite{Li}, Steer-Wren \cite{Steer-Wren}, Mochizuki \cite{Mochizuki} 
and contrasts with the Maruyama-Yokogawa definition which uses a single index \cite{MaruyamaYokogawa}. 

\subsection{Definition:}\label{Defn}
A {\em parabolic bundle} on $(X,D)$ is a collection of vector bundles $F_{\alpha}$  indexed by
multi-indices $\alpha =(\alpha _1,\ldots , \alpha _k)$ with $\alpha _i\in \Q$, together with
inclusions of sheaves of $\cO _X$-modules
$$
F_{\alpha} \hookrightarrow F_{\beta}
$$
whenever $\alpha _i\leq \beta _i$ (a condition which we write as $\alpha \leq \beta$ in what follows),
subject to the following hypotheses:
\newline
---(normalization/support)\, let $\delta ^i$ denote the multiindex
$\delta ^i_i=1,\;\; \delta ^i_j=0, \; i\neq j$,
then $F_{\alpha + \delta ^i} = F_{\alpha } (D_i)$ (compatibly with the inclusion); and
\newline
---(semicontinuity)\, for any given $\alpha$ there exists $c>0$ such that for any multiindex $\varepsilon$
with $0\leq \varepsilon _i < c$ we have $F_{\alpha + \varepsilon} = F_{\alpha}$.

It follows from the normalization/support condition that the quotient sheaves $F_{\alpha} / F_{\beta} $
for $\beta \leq \alpha$  are supported in a schematic neighborhood of the divisor $D$, and indeed if
$\beta \leq \alpha \leq \beta + \sum n_i\delta ^i$ then $F_{\alpha} / F_{\beta}$ is supported over
the scheme $\sum _{i=1}^k n_iD_i$. Let $\delta  := \sum _{i=1}^k\delta ^i$. Then
$$
F_{\alpha - \delta} = F_{\alpha} (-D)
$$
and $F_{\alpha} / F_{\alpha -\delta} = F_{\alpha} |_D$.

The semicontinuity condition means that the structure is determined by the sheaves $F_{\alpha}$ for
a finite collection of indices $\alpha$ with $0\leq \alpha _i < 1$, the {\em weights}.

A parabolic bundle is called \textit{locally abelian} if in a Zariski neighbourhood of any point $x\in X$ 
there is an isomorphism between $F$ and a direct sum of parabolic line bundles.
By Lemma 3.3 of \cite{Iy-Si}, it is equivalent to require this condition on an etale neighborhood. 

The locally abelian condition first appeared in Mochizuki's paper \cite{MochizukiA}, in the form of his notion of {\em compatible filtrations}.
The condition that there be a global frame splitting all of the parabolic filtrations appears as the conclusion of his Corollary 4.4 in \cite[\S 4.5.3]{MochizukiA},
cf Theorem \ref{locabcriterion} below. A somewhat similar compatibility condition appeared earlier in Li's paper \cite[Definition 2.1(a)]{Li},
however his condition is considerably stronger than that of \cite{MochizukiA} and some locally abelian cases such as
Case B in \S \ref{example} below will not be covered by  \cite{Li}. The notion of existence of a local frame splitting all of the
filtrations, which is our definition of ``locally abelian'', did occur as the conclusion of \cite[Lemma 3.2]{Li}. 

Fix a single $n$ which will be the
denominator for all of the divisor components, to make notation easier. Let $m$ be the
number of divisor components, and introduce the notation
$$
[a_1,\ldots , a_m] := (\frac{a_1}{n}, \ldots , \frac{a_m}{n})
$$
for multi-indices, so the parabolic structure is determined by the bundles
$F_{[a_1,\ldots , a_m]}$ for $0\leq a_i < n$ with $a_i$ integers. 


\subsection{Parabolic bundles by filtrations}
\label{se.-parfilt}

Historically the first way of considering parabolic bundles was by filtrations on the restriction to 
divisor components \cite{MehtaSeshadri}, \cite{Seshadri},
see also \cite{MaruyamaYokogawa}, \cite{Biswas}, \cite{InabaEtAl} \cite{Li} \cite{Steer-Wren} \cite{Mochizuki} \cite{Panov}.  
Suppose we have a vector bundle $E$ and filtrations of $E|_{D_i}$ by saturated subbundles:
$$
E|_{D_i}=F^i_{0}\supset F^i_{-1}\supset...\supset F^i_{-n}=0
$$
for each $i$, $1\leq i\leq m$.

Consider the kernel sheaves for $-n \leq j \leq 0$,
$$
0\lrar \ov{F^i_j}\lrar E\lrar \f{E|_{D_i}}{F^i_{j}}\lrar 0
$$
and define
\begin{equation}
\label{intersectionformula}
F_{[a_1,a_2,...,a_m]}:=\cap_{i=1}^m \ov{F^i_{a_i}},
\end{equation}
for $-n \leq a_i \leq 0$. In particular $F_{[0,\ldots , 0]} = E$.
This can then be extended to sheaves defined for all values of $a_i$ using the
normalization/support condition 
\begin{equation}
\label{extensioncond}
F_{[a_1,\ldots , a_i+n,\ldots , a_m]} = F_{[a_1,a_2,...,a_m]}(D_i).
\end{equation}
We call this a parabolic structure {\em given by filtrations}.

Conversely, suppose we are given a parabolic structure $F_{\cdot}$ as described in \eqref{Defn} when all the component sheaves $F_{[a_1,\ldots , a_m]}$ are 
vector bundles. Set $E:= F_{[0,\ldots , 0]}$,
and note that
$$
E|_{D_i} = E / F_{[0,\ldots , -n, \ldots , 0]}
$$
where $-n$ is put in the $i$th place. The image of $F_{[0,\ldots , -a_i,\ldots , 0]}$ in $E|_{D_i}$ is
a subsheaf, and we assume that it is a saturated subbundle. This gives a parabolic structure ``by filtrations''.
We can recover the original parabolic structure $F_{\cdot}$ by the intersection formula \eqref{intersectionformula}.

We feel that these constructions only make good sense under the
locally abelian hypothesis. We note some consequences of the locally abelian property. 

\begin{lemma}
\label{locabconseq}
Suppose $\{F_{[a_1,a_2,...,a_m]}\}_{ -n\leq a_i\leq 0}$ define a locally abelian parabolic bundle on $X$
with respect to $(D_1,\ldots , D_m)$. Let $E:= F_{[0,\ldots ,0]}$, which is a vector bundle
on $X$. Then $F$ comes from a construction as above
using unique filtrations of $E|_{D_i}$ and we have the following properties:
\newline
(a)\,\, the $F_{[a_1,a_2,...,a_m]}$ are locally free;
\newline
(b)\,\, for each $k$ and collection of indices $(i_1,\ldots ,i_k)$, at each point in the $k$-fold intersection
$P\in D_{i_1}\cap \cdots \cap D_{i_k}$ the filtrations $F^{i_1}_{\cdot},\ldots ,F^{i_k}_{\cdot}$ of $E_P$
admit a common splitting, hence the associated-graded 
$$
Gr^{F^{i_1}}_{j_1}\cdots Gr^{F^{i_k}}_{j_k}(E_P)
$$
is independent of the order in which it is taken (see \cite{DeligneHodge2}); and 
\newline
(c)\,\, the
functions 
$$
P\mapsto \dim \, Gr^{F^{i_1}}_{j_1}\cdots Gr^{F^{i_k}}_{j_k}(E_P)
$$
are locally constant functions of $P$ on the multiple intersections $D_{i_1}\cap \cdots \cap D_{i_k}$.
\end{lemma}
\begin{proof}
Direct. 
\end{proof}

The above conditions are essentially what Mochizuki has called ``compatibility'' of the filtrations \cite[\S 4]{MochizukiA}, 
and he shows that
they are sufficient for obtaining a compatible local frame. Compare with \cite[Lemma 3.2]{Li} where the proof is much
shorter because the compatibility condition in the hypothesis is stronger. 

\begin{theorem}[{Mochizuki \cite[\S 4.5.3, Cor. 4.4]{MochizukiA}}]
\label{locabcriterion}
Suppose given a parabolic structure which is a collection of sheaves $F_{[a_1,a_2,...,a_m]}$ obtained from
filtrations on a bundle $E$ as above. If these satisfy
conditions (a), (b) and (c) of the previous lemma, then the parabolic structure is locally abelian. 
\end{theorem}

The situation is simpler in the case of surfaces which we describe here.

\begin{lemma}\label{surface}
Suppose $X$ is a surface with a normal crossings divisor $D=D_1+\ldots + D_m\subset X$. Suppose given data of 
a bundle $E$ and strict filtrations of $E|_{D_i}$ as in Lemma \ref{locabconseq}. Then this data defines a locally abelian parabolic
bundle on $(X,D)$.
\end{lemma}
\begin{proof}
One way to prove this is to use the correspondence with bundles on the DM-stack covering 
$Z:= X\langle \frac{D_1}{n},\ldots ,\frac{D_m}{n}\rangle$ (see \cite[Lemma 2.3]{Iy-Si}).
Let $Z'$ be the complement of the intersection points of the divisor. On $Z'$ the given filtrations
define a vector bundle, as can be seen by applying the 
correspondence of \cite{Borne} \cite{Iy-Si} in codimension $1$, or more
concretely just by using the filtrations to make a sequence of
elementary transformations. Then, since $Z$ is a smooth surface, this bundle extends
to a unique bundle on $Z$, which corresponds to a locally abelian parabolic bundle on $X$ \cite{Iy-Si}. 

Another way to prove this is to note that there are only double intersections. At a point $P$ where $D_i$ and
$D_j$ intersect, the filtrations coming from $D_i$ and $D_j$ have a common splitting. This can then be extended
along both $D_i$ and $D_j$ as a splitting of the respective filtrations, and extended in any way to the rest of $X$.
The resulting direct sum decomposition splits the parabolic structure. This is illustrated by an example in \S \ref{example}.
\end{proof}

We mention here a more general notation used by Mochizuki \cite[\S 3.1]{Mochizuki}
for parabolic bundles given by a filtration,
starting with an origin ${\bf c}= (c_1,\ldots , c_m)$ which may be different from $(0,\ldots , 0)$.
In this case, the underlying bundle is 
$$
E:= F_{[c_1,\ldots , c_m]}
$$
and the filtrations on $E|_{D_i}$ are denoted $F^i_j$ indexed by $c_i-n\leq j \leq c_i$
with $F^i_{c_i}=E|_{D_i}$ and $F^i_{c_i-n}= 0$. We can go between different values of ${\bf c}$ by tensoring
with parabolic line bundles.


\subsection{Parabolic sheaves in the Maruyama-Yokogawa notation}
\label{se.-parmy}

In their original definition of parabolic structures on higher-dimensional varieties,
 Maruyama and Yokogawa considered the general notion of parabolic sheaf with respect to a single
divisor, even if the divisor is not smooth \cite{MaruyamaYokogawa}. Call this a {\em  MY parabolic struture}. 
We can apply their definition to
the full divisor $D = D_1+\ldots + D_m$. This is what was done for example in Biswas 
\cite{Biswas}, Borne \cite{Borne}
and many other places. Of course for the case of curves, the two are completely equivalent because a divisor is always
a disjoint union of its components;
multi-indexed divisors were used by Mehta and Seshadri \cite{MehtaSeshadri}.  
Some of the first places where multi-indexed divisors were used in higher dimensions were in Li \cite{Li}, Steer-Wren \cite{Steer-Wren}, Panov \cite{Panov} and
Mochizuki \cite{Mochizuki}.  In the MY case the parabolic structure is given by a collection of sheaves indexed
by a single parameter $F^{\alpha}$ for $\alpha \in \Q $, with $F^{\alpha + 1} = F^{\alpha} (D)$. We use upper
indexing to distinguish this from our notations (although they would be the same in the case of a single
smooth divisor). If $F_{\cdot}$ is
a parabolic structure according to our notations, then we get a MY-parabolic structure by setting
$$
F^{\alpha} := F_{\alpha , \ldots , \alpha}.
$$
Conversely, given a MY-parabolic structure $F^{\cdot}$, if we assume that $E:=F^0$ is a bundle, then
the images of $F^{-\frac{a_i}{n}}$ in $E_{D_i}$ define subsheaves at generic points of the components $D_i$,
which we can complete to saturated subsheaves everywhere. If $F^{\cdot}$ is locally abelian (that is to say,
locally a direct sum of MY-parabolic line bundles) then these saturated subsheaves are subbundles and we recover
the parabolic structure via filtrations, hence the parabolic structure $F_{\cdot}$ in this way. 
This construction is tacitly used by Biswas in \cite[pp. 599, 602]{Biswas2}, although he formally sticks to the MY-parabolic  notation. 

In the locally
abelian case, all of these different points of view permit us to represent the same objects and going between
them by the various constructions we have outlined, is a commutative process in the sense that by any path
we get back to the
same objects in each notation. We don't attempt to identify the optimal set of hypotheses, weaker than
locally abelian, on the various structures
which would allow to give a more general statement of this sort of commutation of the various constructions.
This doesn't seem immediately relevant since, for now, it doesn't seem clear what is the really good notion
of parabolic sheaf.


\subsection{Parabolic bundles on a DM-stack}

Recall from \cite{Borne} \cite{Cadman} \cite{MatsukiOlsson} \cite{Iy-Si} that given $(X,D)$ and a denominator $n$,
 we can form a
DM-stack denoted $Z:= X\langle  \frac{D_1}{n},\ldots ,  \frac{D_m}{n}\rangle$, and there is an equivalence of categories 
between parabolic bundles on $(X,D)$ with denominator
$n$, and vector bundles on the DM-stack $Z$. The Chern character will be defined using this equivalence, and we would
like to analyse it by an induction on the number of divisor components $m$. Thus, we are interested in intermediate cases
of parabolic bundles on DM-stacks.  

We can carry out all the above constructions in the case when $X$ is a DM stack and $D_i$ are
smooth divisors, i.e., smooth closed substacks of codimension $1$,  meeting transversally on $X$. 

\begin{lemma}\label{le.-stack}
The construction $(X,D)\mapsto Z:= X\langle \frac{D}{n}\rangle $ makes sense for any smooth DM stack $X$ and
smooth divisor $D\subset X$. The stack
$Z$ is then again smooth with a morphism of stacks $Z\rightarrow X$. 
\end{lemma}
\begin{proof}
Since the construction \cite{Cadman} \cite{MatsukiOlsson} \cite{Borne}
of the DM-stack $X\langle \frac{D}{n}\rangle $ when $X$ is a variety 
is local for the \'etale topology (see \cite[\S 2.2]{Iy-Si}), the same construction works when $X$ is a DM-stack.
\end{proof}

Let $Z_k := X\langle \frac{D_1}{n}, \ldots , \frac{D_k}{n}\rangle $. This is a DM-stack (see \cite{Cadman} \cite{MatsukiOlsson} \cite{Borne}
\cite[\S 2.2]{Iy-Si}) and we have maps
$$
\ldots \rightarrow Z_k \rightarrow Z_{k-1} \rightarrow \ldots \rightarrow Z_0 = X.
$$
On $Z_k$ we have divisors $D_j^{(k)}$ which are the pullbacks of the divisors $D_j$ from $X$. 
When $j>k$ the divisor $D_j^{(k)}$ is smooth, whereas for $j\leq k$ the divisor $D_j^{(k)}$ has
multiplicity $n$.

\begin{lemma}
With the above notations, we have the inductive statement that for any $0\leq k < m$, 
$$
Z_{k+1} = Z_k \langle \frac{D_{k+1}^{(k)}}{n}\rangle .
$$
\end{lemma}
\begin{proof}
Recall the definition of $Z_{k+1}$ : if we assume $D_i$ for $i=1,\ldots,k+1$ is defined by equations $z_i=0$ and 
on any local chart (for the \'etale topology) some of the components say $D_1,\ldots, D_{k'}$ occur then
the local chart for $Z_{k+1}$ with coordinates $u_i$ is defined by the equations $z_i=u_i^n$ for $i=1,\ldots,k'$ and $z_i=u_i$ for $i>k'$.
Now $Z_k \langle \frac{D_{k+1}^{(k)}}{n}\rangle $ is obtained from $Z_k$ by defining local chart with coordinates $w_i$
and repeating the above construction by considering the component divisor $D^{(k)}_{k+1}$ on $Z_k$, by applying Lemma \ref{le.-stack} 
and having the same denominator $n$. It is now clear that
both the constructions define the same stack.
\end{proof}

Suppose $X$ is a smooth DM stack and $D\subset X$ is a smooth divisor. Then we define the notion of {\em parabolic
bundle} on $(X,D)$ (with $n$ as denominator) as follows. A parabolic structure 
is a collection of sheaves $F_{\alpha}$ on $X$
(with $\alpha \in \frac{1}{n} \Z$) with $F_{[a]}\rightarrow F_{[a+1]}$ (remember the notation at the start
here with $m=1$ so $[a]= (\frac{a}{n})$).  This is a parabolic bundle if the 
$F_{[a]}$ are bundles
and the quotient sheaves 
$$
F_{[a+1]}/F_{[a]}
$$
are bundles supported on $D$. This is equivalent to a locally abelian 
condition in the \'etale topology
of $X$. Indeed, we can attach weights $\f{a}{n}$ to the graded pieces $F_{[a+1]}/F_{[a]}$ 
whenever this is non-zero and define locally on a general point of the divisor $D$ a direct sum $L$ of parabolic line bundles 
such that if the rank of $F_{[a+1]}/F_{[a]}$ is 
$n_a$ then $L= \sum_a \cO(-\f{a}{n}D)^{\oplus n_a}$.    

\begin{lemma}
There is an equivalence of categories between bundles on $X\langle \frac{D}{n}\rangle $ and parabolic bundles on $(X,D)$
with $n$ as denominator.
\end{lemma}
\begin{proof}
This is proved by Borne \cite[Theorem 5]{Borne} when $X$ is a smooth variety. In the case of a DM stack since everything is local in the \'etale 
topology it works the same way.
\end{proof}

Similarly if $D_i$ are smooth divisors meeting transversally on a DM stack $X$ then we can define a notion
of locally abelian parabolic bundle on $(X; \sum _i D_i)$, as in \S 2.1. Here the locally abelian condition is local in the etale
topology which is the only appropriate topology to work with on $X$. 

\begin{lemma}
\label{equivchain}
With the notations of the beginning, the categories of
locally abelian parabolic bundles on 
$$
(Z_k; D_{k+1}^{(k)}, \ldots , D_m^{(k)})
$$
are all naturally equivalent.
\end{lemma}
\begin{proof}
When $k=m$ and for any $k$, so we consider $Z_m$ and $Z_k$, the equivalence of vector bundles on $Z_m$ and locally abelian parabolic bundles on $Z_k$ is proved in 
\cite[Lemma 2.3]{Iy-Si} (actually it is proved when $Z_k$ is a variety but as earlier the same proof holds for the DM-stack $Z_k$). 
This gives the equivalences of categories on any $Z_k$ and $Z_{k'}$. 
\end{proof}

In particular the cases $k=0$ so $Z_k=X$ and $k=m$ where there are no further divisor components, correspond to the
equivalence of categories of \cite[Lemma 2.3]{Iy-Si}:

\begin{corollary}
\label{equivZm}
The category of locally abelian parabolic bundles on $X$ is
equivalent to the category of vector bundles on 
 $Z_m = X\langle \frac{D_1}{n},\ldots , \frac{D_m}{n}\rangle$.
\end{corollary}

\subsection{Chern characters}

We recall here the abstract definition of the Chern character of a parabolic bundle. If $F_{\cdot}$ is
a parabolic bundle with rational weights having common denominator $n$, then it corresponds to
a vector bundle $F_{DM}$ on the DM-stack $Z_m = X\langle \frac{D_1}{n},\ldots ,\frac{D_m}{n}\rangle $.
Let $\pi : Z_m \rightarrow X$ denote the projection. By Gillet \cite{Gi} and Vistoli \cite{Vistoli} it induces an isomorphism
of rational Chow groups
\begin{equation}
\label{chowidentification}
\pi _{\ast} : CH(Z_m)_{\Q } \stackrel{\cong}{\longrightarrow} CH(X)_{\Q }.
\end{equation}
In \cite{Iy-Si}, following an idea of Borne \cite{Borne}, we defined the {\em Chern character of $F$} to be
\begin{equation}
\label{chdef}
\m{ch}(F):= \pi _{\ast}(\m{ch} (F_{DM})) \in CH(X)_{\Q }.
\end{equation}
It is a formal consequence of this definition that Chern character is compatible additively with direct sums (or more
generally extensions), multiplicatively with tensor products, and the pullback of the Chern character is the Chern
character of the pullback bundle for a morphism $f$ of varieties if the normal-crossings divisors are in standard position
with respect to $f$.


\subsection{Statement of the main theorem}

Our goal is to give a formula for the Chern character defined abstractly by \eqref{chdef}.
The first main theorem is that the Chern character depends only on the Chern characters of
the component bundles, and not on the inclusion morphisms between them. This is not in any way
tautological, as is shown by the examples we shall consider in \S \ref{formula} below which show that
it is not enough to consider the Chern characters of the bundle $E$ plus the filtrations,
or just the Maruyama-Yokogawa components. The full collection of component
bundles $F_{[a_1,\ldots , a_m]}$ is sufficient to account for the incidence data among the filtrations,
and allows us to obtain the Chern character. 

\begin{theorem}
\label{invariance-th}
Suppose $F$ and $G$ are locally abelian parabolic bundles on a DM stack $X$ with $n$ as denominator.
Suppose that for all $a_i$ with $0\leq a_i <n$ the bundles 
$F_{[a_1,\ldots , a_m]}$  and $G_{[a_1,\ldots , a_m]}$ have the same Chern character in the rational
Chow groups of $X$. Then the parabolic bundles
$F$ and $G$ have the same Chern character in the rational Chow group of $X$.
\end{theorem}

When we have two parabolic bundles $F$ and $G$ satisfying the hypothesis of the theorem, we say that
$F$ and $G$ are {\em componentwise Chow equivalent}. A stronger condition is 
to say that $F$ and $G$ are {\em componentwise isomorphic}, meaning that the 
$F_{[a_1,\ldots , a_n]}$  and $G_{[a_1,\ldots , a_n]}$ are isomorphic bundles on $X$. This obviously
implies that they are componentwise Chow equivalent, and so the theorem will imply that they have the 
same Chern character. 

Once we have Theorem \ref{invariance-th}, it is relatively straightforward to give an explicit calculation
of the Chern character by exhibiting a componentwise isomorphism of parabolic bundles. The componentwise
isomorphism which will come into play, will not, however, come from an isomorphism of parabolic structures 
because the individual
isomorphisms on component bundles will not respect the inclusion maps in the parabolic 
structure. The resulting formula is a weighted average as stated in Theorem \ref{mainformula},
proven as Theorem \ref{answer} below.


\section{Reduction to the case of one divisor}


In this section and the next, we prove Theorem \ref{invariance-th}.
In this section we will use the intermediate stacks
$Z_k$ in order to reduce to the case of only one smooth divisor component; then in the next section we prove the formula
for that case. To see how the reduction works we
have to note what happens to the component bundles in the equivalence of Lemma \ref{equivchain}.

Fix $0<k\leq m$ and consider the equivalence of Lemma \ref{equivchain} which we denote (a) in what follows: suppose $E$ is a 
locally abelian parabolic bundle on
$$
(Z_{k-1}; D_{k}^{(k-1)},\ldots , D_m^{(k-1)}),
$$
then it corresponds to $F$ which is a locally abelian parabolic bundle on 
$$
(Z_k; D_{k+1}^{(k)},\ldots , D_m^{(k)}).
$$
Recall that $Z_k = Z_{k-1}\langle \frac{D_k^{(k-1)}}{n}\rangle $ and that we have an equivalence (b) between bundles
on $Z_k$, and parabolic bundles on $Z_{k-1}$ with respect to the divisor $D_k^{(k-1)}$.
For any $b_{k+1},\ldots , b_m$ we can let $a_k$ vary, and using $E$ we obtain a parabolic bundle
$$
H^{[b_{k+1},\ldots , b_m]}:= a_k \mapsto E_{[a_k,b_{k+1},\ldots , b_m]}
$$
on $Z_{k-1}$ with respect to the divisor $D_k^{(k-1)}$. 

\begin{lemma}
\label{components}
Suppose that $E$ and $F$ correspond via the equivalence (a) as in the above notations, and define the
parabolic bundle $H^{[b_{k+1},\ldots , b_m]}$ as above, which for any $b_{k+1},\ldots , b_m$ is a parabolic
bundle on $Z_{k-1}$ with respect to the divisor $D_k^{(k-1)}$.  Then this parabolic bundle 
$H^{[b_{k+1},\ldots , b_m]}$ is the one which corresponds via the equivalence (b) to the component vector bundle
$F_{[b_{k+1},\ldots , b_m]}$ of the parabolic bundle $F$.
\end{lemma}
\begin{proof}
We use the definition of the pushforward (\cite[\S 2.2]{Iy-Si}) which provides the explicit equivalence
in Lemma \ref{equivchain}. For simplicity, we assume that $k=m-1$ so we are looking at the case $$Z_m\sta{p}{\lrar} Z_{m-1}\sta{q}{\lrar} Z_{m-2}.$$
Let $G$ be the vector bundle on $Z_m$ corresponding to $E$ or $F$, using the equivalence in Lemma \ref{equivchain}. 
Consider the vector bundle $F_{[b_m]}$ on $Z_{m-1}$. We want to check that the associated parabolic bundle 
$q_*F_{[b_m]}$ is $H^{b_m}$. The following equalities prove this claim.
\begin{eqnarray*}
(q_*F_{[b_m]})_{a_{m-1}} &=& q_*(F_{[b_m]}(a_{m-1}D_{m-1})) \\ 
                         &=& q_*((p_*G)_{[b_m]}(a_{m-1}D_{m-1})) \\
                         &=& q_*(p_*(G(b_mD_m)(a_{m-1}D_{m-1}))\\
                         &=& (q\circ p)_*(G(b_mD_m +a_{m-1}D_{m-1}))\\
                         &=& E_{[a_{m-1},b_m]}\\
                         &=& (H^{[b_m]})_{[a_{m-1}]}
\end{eqnarray*}

\end{proof}

A corollary of this observation is that we can reduce for Theorem \ref{invariance-th} to the case of a
single divisor. 

\begin{corollary}
\label{reduction}
Suppose that Theorem \ref{invariance-th} is known for $m=1$, that is for a single smooth divisor. 
Then it holds in general.
\end{corollary}
\begin{proof}
Fix $X$ with $D_1,\ldots , D_m$ and define the sequence of intermediate stacks $Z_k$ as above.
Suppose $F$ and $G$ are locally abelian parabolic bundles on $X=Z_0$ which are componentwise Chow 
equivalent. For any $k$ let $F^{(k)}$ and $G^{(k)}$ denote the corresponding locally abelian 
parabolic bundles on $Z_k$ with respect to the remaining divisors $D_{k+1}^{(k)},\ldots , D_m^{(k)}$.
We claim by induction on $0\leq k \leq m$ that the $F^{(k)}$ and $G^{(k)}$ are componentwise Chow equivalent.
This is tautologically true for $k=0$. Fix $0<k\leq m$ and suppose it is true for $k-1$.
Then $F^{(k-1)}$ and $G^{(k-1)}$ induce for any $b_{k+1},\ldots , b_m$ parabolic bundles which
we can denote by $H_F^{[b_{k+1},\ldots , b_m]}$ and $H_G^{[b_{k+1},\ldots , b_m]}$, as in 
Lemma \ref{components}.
These are parabolic bundles on $Z_{k-1}$ with respect to the single smooth divisor $D_k^{(k-1)}$.  
The components of these parabolic bundles are Chow equivalent, since they come from the components
of $F^{(k-1)}$ and $G^{(k-1)}$ which by the induction hypothesis are componentwise Chow equivalent.
Therefore, considered as parabolic bundles with respect to a single divisor, 
$H_F^{[b_{k+1},\ldots , b_m]}$ and $H_G^{[b_{k+1},\ldots , b_m]}$ are componentwise Chow equivalent.
In the present corollary we are assuming that Theorem \ref{invariance-th} is known for the case 
$m=1$ of a single divisor. Applying this case of Theorem \ref{invariance-th} we get that the bundles on $Z_k$ associated to 
$H_F^{[b_{k+1},\ldots , b_m]}$ and $H_G^{[b_{k+1},\ldots , b_m]}$ are Chow equivalent. 
However, by Lemma \ref{components} applied to the comparison between 
$F^{(k-1)}$ and $F^{(k)}$, the bundle on $Z_k$ corresponding to the
parabolic bundle $H_F^{[b_{k+1},\ldots , b_m]}$ is exactly the component
$$
F^{(k)}_{[b_{k+1},\ldots , b_m]}.
$$
Similarly, applying Lemma \ref{components} to the comparison between 
$G^{(k-1)}$ and $G^{(k)}$, the bundle on $Z_k$ corresponding 
to the
parabolic bundle $H_G^{[b_{k+1},\ldots , b_m]}$ is exactly the component
$$
G^{(k)}_{[b_{k+1},\ldots , b_m]}.
$$
Thus the result of our application of the single divisor case of Theorem \ref{invariance-th} 
is that the bundles $F^{(k)}_{[b_{k+1},\ldots , b_m]}$ and 
$G^{(k)}_{[b_{k+1},\ldots , b_m]}$ are Chow equivalent. This exactly says that the parabolic bundles
$F^{(k)}$ and $G^{(k)}$ are componentwise Chow equivalent, which completes our induction step.
\end{proof}


\section{The single divisor case}


By Corollary \ref{reduction}, it now suffices to prove Theorem \ref{invariance-th} in the case $m=1$. 
Simplify notation. Suppose we have a smooth DM stack $X$ and a smooth divisor $D$,
 and suppose we have a parabolic bundle $F$ on $X$ with respect to $D$.
It is a collection of bundles denoted $F_{[a]}$ with $a\in \Z$ (as usual without saying so we
assume that the denominator is $n$). Let $Z:= X\langle \frac{D}{n}\rangle $, so $F$ 
corresponds to a vector bundle $E$ on $Z$. According to the definition \eqref{chdef} we would like to
show that the Chern character of $E$ in the rational Chow group of $Z$ depends only on the Chern characters
of the $F_{[a]}$ in the rational Chow group of $X$, noting the identification \eqref{chowidentification}.

Let $p:Z\rightarrow X$ denote the map of DM stacks. The inverse image $p^{\ast}(D)$ is a divisor 
in $Z$ which has multiplicity $n$, because $p$ is totally ramified of degree $n$ over $D$.
In particular, there is a divisor $R\subset Z$ such that 
$$
p^{\ast}(D) = n\cdot R.
$$
This $R$ is well-defined as a smooth closed substack of codimension $1$ in $Z$. However, $R$ is a gerb over $D$.
More precisely, we have a map $R\rightarrow D$ and there is a covering of $D$ in the etale topology 
by maps $U\rightarrow D$ such that there is a lifting $U\rightarrow R$. If we are given such a lifting then
this gives a trivialization 
$$
U\times _D R \cong U\times B(\Z /n),
$$
where $B(\Z /n)$ is the one-point stack with group $\Z /n$. This can be summed up by saying that $R$ is a gerb
over $D$ with group $\Z /n$. It is in general not trivial. (We conjecture that the obstruction is the same
as the obstruction to the normal bundle $N_{D/X}$ having an $n$-th root as line bundle on $D$.)
On the other hand, the character theory for $R$ over $D$ is trivialized in the following sense. There is a 
line bundle $N := \cO _X (R)|_R$ on $R$ with the property that on any fiber of the form $B(\Z /n)$, 
$N$ is the primitive character of $\Z /n$.  

Using $N$, we get a canonical decomposition of bundles on $R$. Suppose $E$ is a bundle on $R$. Then
$p_{R,\ast}E$ is a bundle on $D$ which corresponds in each fiber to the trivial character. Here
$p_R$ is the map $p$ restricted to $R$. For any $i$
we have a map 
$$
p_R^{\ast}(p_{R,\ast}(E\otimes N^{\otimes -i})) \otimes N^{\otimes i} \rightarrow E.
$$

\begin{lemma}
\label{decompR}
If $E$ is a bundle on $R$ then the above maps put together for $0\leq i < n$ give a direct sum decomposition
$$
\bigoplus _{i=0}^{n-1} p_R^{\ast}(p_{R,\ast}(E\otimes N^{\otimes -i})) \otimes N^{\otimes i} 
\stackrel{\cong}{\longrightarrow} E.
$$
\end{lemma}
\begin{proof}
The maps exist globally. To check that the map is an isomorphism it suffices to do it locally over $D$ in
the etale topology (since the map $p_R$ is involved). As noted above, locally over $D$ the gerb $R$
is a product of the form $U\times B(\Z /n)$. A bundle $E$ on the product is the same thing as a bundle on $U$ 
together with an action of the group $\Z / n$. In turn this is the same thing as a bundle with action
of the group algebra $\cO _U [\Z / n]$ but relative $Spec$ of this algebra over $U$ is a disjoint union of
$n$ copies of $U$, so $E$ decomposes as a direct sum of pieces corresponding to these sections. 
This decomposition may be written as $E = \bigoplus _{\chi} E_{\chi}$ where the $\chi$ are characters of $\Z / n$
and $\Z / n$ acts on $E_{\chi}$ via the character $\chi$. In terms of the DM stack this means that $E$ decomposes
as a direct sum of bundles on $U$
tensored with characters of $\Z /n$ considered as line bundles on $B(\Z /n)$. 
Using this decomposition
we can check that the above map is an isomorphism (actually it gives back the same decomposition). 
\end{proof}

Now suppose $E$ is a bundle on $Z$. Then its restriction to $R$, noted $E_R$, decomposes according to the above
lemma. Define two pieces as follows: $E_{R,{\rm fix}}$ is the piece corresponding to $i=0$ in the decomposition.
Thus
$$
E_{R,{\rm fix}} = p_R^{\ast}(p_{R,\ast}E_R).
$$
On the other hand, let $E_{R,{\rm var}}$ denote the direct sum of the other pieces in the decomposition.
The decomposition of Lemma \ref{decompR} thus gives a direct sum decomposition
$$
E_R =  E_{R,{\rm fix}}\oplus E_{R,{\rm var}}.
$$

Define the {\em standard elementary transformation} ${\bf e}(E)$ of a bundle $E$ over $Z$, as the kernel
\begin{equation}
\label{elementtransform}
0 \rightarrow {\bf e}(E) \rightarrow E \rightarrow E_{R,{\rm var}} \rightarrow 0.
\end{equation}

\begin{lemma}
\label{reseE}
Suppose $E$ is a bundle on $Z$. Then we have the following exact sequence for the restriction of
the standard elementary transformation of $E$: 
$$
0 \rightarrow E_{R,{\rm var}} \otimes N^{\ast} \rightarrow ({\bf e}(E))_R \rightarrow E_{R,{\rm fix}} \rightarrow 0.
$$
\end{lemma}
\begin{proof}
Consider the exact sequence :

$$0\lrar E\otimes \cO(-R)\lrar E \lrar E_R\lrar 0.$$
 
Since $E_R=E_{R,{\rm fix}}\oplus E_{R,{\rm var}}$, and ${\bf e}(E)$ is the kernel of the composed map $$E\lrar E_R\lrar E_{R,{\rm var}}$$ there is an induced injective
map $$E\otimes \cO(-R)\lrar {\bf e}(E)$$ inducing the restriction map on $R$
$$(E_{R,{\rm fix}}\oplus E_{R,{\rm var}})\otimes \cO(-R)_{|R}\lrar {\bf e}(E)_{|R}$$
The kernel of the restriction
$$({\bf e}(E))_R\lrar E_R\lrar E_{R,{\rm fix}}$$
is clearly $E_{R,{\rm var}} \otimes \cO(-R)_{|R}\, =\,  E_{R,{\rm var}} \otimes N^{\ast}$.

\end{proof}

Suppose $E$ is a bundle on $Z$. Define $\rho (E)$ to be the largest integer $k$ with $0\leq k < n$ such that the piece
$$
p_R^{\ast}(p_{R,\ast}(E_R\otimes N^{\otimes -k})) \otimes N^{\otimes k} 
$$
in the decomposition of Lemma \ref{decompR} is nonzero. 

Actually we may consider this definition for any vector bundle on $R$.

\begin{corollary}
\label{decrease}
The invariant $\rho$ decreases under the standard elementary transformation: if $\rho (E) > 0$ then
$$
\rho ({\bf e}(E)) < \rho (E).
$$
\end{corollary}
\begin{proof}
Consider the exact sequence from Lemma \ref{reseE} :
$$
0 \rightarrow E_{R,{\rm var}} \otimes N^{\ast} \rightarrow ({\bf e}(E))_R \rightarrow E_{R,{\rm fix}} \rightarrow 0.
$$
Using the pushforward and pullback operations on this exact sequence,
after twisting by powers of $N$, we notice that it suffices to check that $\rho(E_{R,\rm{var}}\otimes N^*)< \rho(E)$ and $\rho(E_{R,\rm{fix}})=0$.
 
Now 
\begin{eqnarray*}
p_R^*p_{R\, *}(E_{R,\rm{fix}}\otimes N^{-k})\otimes N^{k} &=&p_R^*p_{R\,*}(p_R^*p_{R\,*}E\otimes N^{-k})\otimes N^{k} \\
&=& p_R^*(p_{R\,*}E \otimes p_{R\,*}N^{-k})\otimes N^{k} \\
&=& 0 \,\,\rm{ if }\,\,k\neq 0.
\end{eqnarray*}

Also,
\begin{eqnarray*}
& p_R^*p_{R\,*}(E_{R, \rm{var}}\otimes N^{-1}\otimes N^{-k})\otimes N^{k} \\
= & p_R^*p_{R\,*}\left((\sum_{i=1}^{\rho(E)} p_R^*p_{R\,*}(E_R\otimes N^{-i})\otimes N^{i})\otimes N^{-1}\otimes N^{-k}\right)\otimes N^{k} \\
= & p_R^*\left(\sum_{i=1}^{\rho(E)}p_{R\,*}(E_R\otimes N^{-i})\otimes p_{R\,*}N^{i-1-k}\right) \otimes N^{k}. \\
\end{eqnarray*}

The summands in the above term corresponding to $i-1-k\neq 0$ are zero.
 In other words, the only term left is for $i=k+1$, but if $k \geq \rho(E)$ then this doesn't occur and the whole 
is zero. Hence $\rho(E_{R,\rm{var}}\otimes N^*)< \rho(E)$.
\end{proof}

We now describe the pieces in the decomposition of Lemma \ref{decompR} for $E_R$ in terms of the parabolic structure
on $X$. Introduce the following notation: if $F$ is a parabolic bundle on $X$ along the divisor $D$,
then for any $a\in \Z$ set ${\bf gr}_{[a]}(F):= F_{[a]}/F_{[a-1]}$. It is a vector bundle on the divisor $D$.

\begin{lemma}
\label{description}
Suppose $E$ is a bundle on $Z$ corresponding to a parabolic bundle $F$ over $X$. Then for any $a\in \Z$ we have
$$
p_{R,\ast}(E_R \otimes N^{\otimes a}) \cong {\bf gr}_{[a]}(F).
$$
\end{lemma}
\begin{proof}
We have 
$$
F_{[a]} = p_{\ast}(E(aR)).
$$
Note that $R^1p_{\ast}$ vanishes on coherent sheaves, since $p$ is a finite map in the etale topology.
Thus $p_{\ast}$ is exact. This gives
$$
{\bf gr}_{[a]}(F) = p_{\ast}(E\otimes (\cO _Z(aR) / \cO _Z((a-1)R))).
$$
However, $(\cO _Z(aR) / \cO _Z((a-1)R))$ is a bundle on $R$ which is equal to $N^{\otimes a}$. 
This gives the statement.

\end{proof}

We say that two bundles on $R$ are {\em Chow equivalent relative to $Z$} if their Chern characters map to 
the same thing in the rational Chow group of $Z$. Caution: this is different from their being Chow equivalent
on $R$, because the map $CH(R)_{\Q}\rightarrow CH(Z)_{\Q}$ might not be injective.

\begin{lemma}\label{pushforward}
Suppose $p:Z=X\langle \f{1}{n}\rangle \lrar X$ is a morphism of DM-stacks as in the beginning of this section. Then the following diagram commutes :
\begin{eqnarray*}
CH^{\cdot}(R)_{\Q} &\lrar & CH^{\cdot}(Z)_{\Q} \\
\downarrow \scriptstyle{\cong} & & \downarrow {\scriptstyle \cong} \\
CH^{\cdot}(D)_{\Q}&\lrar & CH^{\cdot}(X)_{\Q}
\end{eqnarray*}
\end{lemma}
\begin{proof}
Use composition of proper pushforwards \cite{Voisin}. The vertical isomorphisms come from the fact that $Z\rightarrow X$ and $R\rightarrow D$
induce isomorphisms of coarse moduli schemes, and \cite{Vistoli} \cite{Gi}.
\end{proof}

\begin{corollary}
\label{chowinv1}
Suppose $E$ and $G$ are vector bundles on $Z$ corresponding to parabolic bundles $F$ and $H$ respectively on $X$.
If $F$ and $H$ are componentwise Chow equivalent then each of the components in the decompositions of
Lemma \ref{decompR} for $E_R$ and $G_R$ are Chow equivalent relative to $Z$. 
\end{corollary}
\begin{proof}
Since $F$ and $H$ are componentwise Chow equivalent the graded pieces ${\bf gr}_{[a]}(F)$ and ${\bf gr}_{[a]}(H)$ are Chow equivalent on $X$.
Hence by Lemma \ref{description}, $p_{R*}(E\otimes N^a)$ and $p_{R*}(G\otimes N^a)$ are Chow equivalent on $X$,
in other words they  are vector bundles on $D$ 
which are Chow equivalent relative to $X$. The pullback of Chow equivalent
objects on $D$ relative to $X$ are Chow equivalent objects on $R$ 
relative to $Z$, by Lemma \ref{pushforward}. Thus, in the sum decomposition of $E_R$ and $G_R$ as in 
Lemma \ref{decompR},
we conclude that the component sheaves are Chow equivalent relative 
to $Z$.
\end{proof}

\begin{corollary}
\label{chowinv2}
Suppose $E$ and $G$ are vector bundles on $Z$ corresponding to parabolic bundles $F$ and $H$ respectively on $X$.
Suppose that $F$ and $H$ are componentwise Chow equivalent. Then the
sheaves $E_{R,{\rm fix}}$ and $G_{R,{\rm fix}}$ are Chow equivalent on $Z$. Similarly, the sheaves
$E_{R,{\rm var}}$ and $G_{R,{\rm var}}$ are Chow equivalent on $Z$. 
\end{corollary}
\begin{proof}
These sheaves come from the components of the decomposition for $E_R$ and $G_R$.
\end{proof}

\begin{lemma}
\label{chowinv3}
Suppose $E$ and $G$ are vector bundles on $Z$ corresponding to parabolic bundles $F$ and $H$ respectively on $X$.
Suppose that $F$ and $H$ are componentwise Chow equivalent. As a matter of notation, let 
${\bf e}_XF$ and ${\bf e}_XH$ denote the parabolic bundles on $X$ corresponding to the vector bundles
${\bf e}E$ and ${\bf e}G$. Then ${\bf e}_XF$ and ${\bf e}_XH$ are componentwise Chow equivalent. 
\end{lemma}
\begin{proof}
Firstly, we claim that
\begin{equation}\label{claim}
({\bf e}_XF)_{[0]} = F_{[0]}.
\end{equation}
To prove the claim, note that $F_{[0]} = p_{\ast}(E)$. On the other hand, since $E_{R,{\rm var}}$ has only
components which have trivial direct images, we have $p_{\ast}(E_{R,{\rm var}})=0$, so the left exact
sequence for the direct image of \eqref{elementtransform}, shows that 
$$
p_{\ast}({\bf e}E) = p_{\ast}(E).
$$
This gives the claim. 

The same claim holds for $H$. 

Now twist the exact sequence in Lemma \ref{reseE} by $N^a$, 
and take the pushforward (which is exact). Do this for both bundles $E$ and $G$, yielding the exact sequences
$$
0\lrar  p_{R *}(E_{R,\rm{var}}\otimes N^{-1+a}) \lrar p_*({\bf e}(E)_R\otimes N^a) \lrar p_*( E_{R,{\rm fix}} \otimes N^a) \lrar 0
$$
and
$$
0\lrar  p_{R *}(G_{R,\rm{var}}\otimes N^{-1+a}) \lrar p_*({\bf e}(G)_R\otimes N^a) \lrar p_*(G_{R,{\rm fix}} \otimes N^a) \lrar 0.
$$
By the hypothesis, Corollary \ref{chowinv1} applies to say that the various components in the decomposition of Lemma \ref{decompR} for
$E_{R,\rm{var}}$ and $E_{R,{\rm fix}}$ are Chow equivalent relative $Z$ to the corresponding components of
$G_{R,\rm{var}}$ and $G_{R,{\rm fix}}$. Thus the left and right terms of both exact sequences are Chow equivalent relative to $X$,
so $p_*({\bf e}(E)_R\otimes N^a)$ and $p_*({\bf e}(G)_R\otimes N^a)$ are Chow equivalent relative to $X$. 

Hence by Lemma \ref{description}, ${\bf gr}_{[a]}({\bf e}_XF)$ and ${\bf gr}_{[a]}({\bf e}_XH)$ are Chow equivalent relative to $X$.
Together with the above claim \eqref{claim}, we deduce that the constituent bundles of ${\bf e}_XF$ and ${\bf e}_XH$ are Chow equivalent on $X$.
\end{proof}

We can iterate the operation of doing the elementary transform, denoted $E\mapsto {\bf e}^pE$.
This corresponds to a parabolic bundle on $X$ denoted by $F\mapsto {\bf e}_X^pF$. Note that 
this is indeed the iteration of the notation ${\bf e}_X$

\begin{exercise}
Give an
explicit description of ${\bf e}_X$ in terms of parabolic bundles.
\end{exercise}

Because the invariant $\rho (E)$ decreases under the operation of doing the standard
elementary transform (until we get to $\rho =0$) it follows that $\rho ({\bf e}^pE)=0$ for some
$p\leq n$. 

\begin{lemma}
\label{rhozero}
Suppose $E$ is a bundle on $E$ with $\rho (E)=0$. Then $E$ is the pullback of a bundle on $X$.
\end{lemma}
\begin{proof}
In this case, we have $E_R\simeq p_R^*p_{R *}E_R$. Hence by Lemma \ref{description}, it follows that ${\bf gr}_{[a]}(F)=0$ for $a>0$. This implies that
$F$ has only one constituent bundle $F_{[0]}$ and is a usual bundle on $X$.
 Hence $E$ is the pullback of $F_{[0]}$.

\end{proof}

The next lemma gives the induction step for the proof of the theorem.

\begin{lemma}
\label{step}
Suppose $E$ and $G$ are vector bundles on $Z$ corresponding to parabolic bundles $F$ and $H$ respectively on $X$.
Suppose that $F$ and $H$ are componentwise Chow equivalent. Suppose also that
${\bf e}E$ and ${\bf e}G$ are Chow equivalent on $Z$. Then $E$ and $G$ are Chow equivalent on $Z$. 
\end{lemma}
\begin{proof}
The componentwise Chow equivalence gives from Corollary \ref{chowinv2}
that $E_{R,{\rm var}}$ and $G_{R,{\rm var}}$ are Chow equivalent relative
to $Z$. The exact sequence of Lemma \ref{reseE} gives that $E$ and $G$ are Chow equivalent on $Z$. 
\end{proof}

Finally we can prove Theorem \ref{invariance-th} in the single divisor case.

\begin{theorem}
\label{invariance1}
Suppose $E$ and $G$ are vector bundles on $Z$ corresponding to parabolic bundles $F$ and $H$ respectively on $X$.
Suppose that $F$ and $H$ are componentwise Chow equivalent. Then $E$ and $G$ are Chow equivalent on $Z$.
\end{theorem}
\begin{proof}
Do this by descending induction with respect to the number $p$ given above Lemma \ref{rhozero}. 
There is some $p_0$ such that $\rho ({\bf e}^{p_0}E)=0$ and $\rho ({\bf e}^{p_0}G)=0$. These come from bundles
on $X$. By Lemma \ref{chowinv3}, these bundles (which are the zero components of the corresponding parabolic
bundles) are Chow equivalent. Thus ${\bf e}^{p_0}E$ and ${\bf e}^{p_0}G$ are Chow equivalent.
On the other hand, by Lemma \ref{chowinv3}, all of the ${\bf e}^{p}_XF$ and ${\bf e}^{p}_XH$
are componentwise Chow equivalent. It follows from Lemma \ref{step}, if we know that 
${\bf e}^{p+1}E$ and ${\bf e}^{p+1}G$ are Chow equivalent then we get that 
${\bf e}^{p}E$ and ${\bf e}^{p}G$ are Chow equivalent. By descending induction on $p$ we get that $E$ and $G$ 
are Chow equivalent. 
\end{proof}

Using Corollary \ref{reduction}, we have now completed the proof of Theorem \ref{invariance-th}.


\section{A formula for the parabolic Chern character}\label{formula}


Now we would like to use Theorem \ref{invariance-th} to help get a formula for the Chern classes. 
Go back to the general situation of a smooth variety $X$ with smooth divisors $D_1,\ldots , D_m$ intersecting
transversally. 
Once we know the formula for the Chern character of a line bundle, we will no longer need to use the stack $Z=
X[\frac{D_1}{n},\ldots , \frac{D_m}{n}]$.

\begin{lemma}
\label{formulae}
Let $F$ be a parabolic bundle on $X$ with respect to the divisors $D_1,\ldots , D_m$.
Then we can form the twisted parabolic bundle $F\otimes \cO (\sum _{i=0}^{m}\frac{b_i}{n}D_i)$.
We have the formulae
$$
\left( F\otimes \cO (\sum _{i=0}^{m}\frac{b_i}{n}D_i)\right) _{[a_1,\ldots , a_m]} = 
F_{[a_1+b_1,\ldots , a_m+b_m]}
$$
and 
$$
{\rm ch} \left( F\otimes \cO (\sum _{i=0}^{m}\frac{b_i}{n}D_i)\right) = 
e^{\sum _{i=0}^{m}\frac{b_i}{n}D_i} {\rm ch}(F).
$$
\end{lemma}
\begin{proof}
Consider the projection 
$p:Z=X\langle \frac{D_1}{n},\ldots , \frac{D_m}{n}\rangle \lrar X$. Let $E$ be the vector bundle on $Z$ corresponding 
to $F$ on $X$ and $\cO(\sum_ib_iR_i)$ be the line bundle on $Z$ corresponding to $\cO(\sum_i\f{b_i}{n}D_i)$ on $X$. Here
$R_i$ denotes the divisor on $Z$ such that $p^*D_i=n.R_i$.

Notice that
\begin{eqnarray*}
\left(F\otimes \cO(\sum_i \f{b_i}{n}D_i)\right)_{[a_1,\ldots,a_m]}&=& p_*(E\otimes \cO(\sum_ib_iR_i)\otimes \cO(a_iR_i))\\
&=& p_*(E\otimes \cO(\sum_i (a_i+b_i) R_i))\\
&=& F_{[a_1+b_1,\ldots , a_m+b_m]}.
\end{eqnarray*} 
The formula for the Chern character is due to the fact that the Chern character defined as we are doing through DM-stacks
is multiplicative for tensor products,
and coincides with the exponential for rational divisors, see \cite{Iy-Si}. 
\end{proof}

\begin{lemma}
\label{trivial}
We have the formula for the trivial line bundle $\cO$ considered as a parabolic bundle:
$$
\cO _{[a_1,\ldots , a_m]} = \cO (\sum _{i=0}^m [\frac{a_i}{n}] D_i)
$$
where the square brackets on the right signify the greatest integer function (on the left they are the
notation we introduced at the beginning). 
\end{lemma}
\begin{proof}
This follows from the definition as in Lemma \ref{formulae}.
\end{proof}

\begin{corollary}
\label{twistcomponents}
Suppose $E$ is a vector bundle on $X$ considered as a parabolic bundle with its trivial structure. Then
$$
\left( E\otimes \cO (\sum _{i=0}^{m}\frac{b_i}{n}D_i)\right) _{[a_1,\ldots , a_m]}
= E (\sum _{i=0}^m [\frac{a_i+b_i}{n}] D_i).
$$
\end{corollary}
\begin{proof}
Use the definition of associated parabolic bundle as in Lemma \ref{formulae}.
\end{proof}

Suppose $F$ is a parabolic bundle on $X$ with respect to $D_1,\ldots , D_m$. 
We will now show by calculation that the two parabolic bundles
$$
\left( \bigoplus _{k_1=0}^{n-1}\cdots \bigoplus _{k_m=0}^{n-1} \cO (-\sum _{i=1}^m \f{k_i}{n}D_i)\right) \otimes F
$$
and
$$
\left( \bigoplus _{u_1=0}^{n-1}\cdots \bigoplus _{u_m=0}^{n-1} F_{[u_1,\ldots , u_m]}\otimes 
\cO (-\sum _{i=1}^m \f{u_i}{n}D_i)\right) 
$$
are componentwise isomorphic (and hence, componentwise Chow equivalent). Notice that the second bundle is
a direct sum of vector bundles on $X$, the component bundles of $F$, tensored with parabolic line bundles, whereas
the first is $F$ tensored with a bundle of positive rank. This will then allow us to get a formula for 
${\rm ch}(F)$.

\begin{lemma}
\label{firstcomponents}
For any $0\leq a_i <n$ we have
$$
\left(
\left( \bigoplus _{k_1=0}^{n-1}\cdots \bigoplus _{k_m=0}^{n-1} \cO (-\sum _{i=1}^m \f{k_i}{n}D_i)\right) \otimes F \right)_{[a_1,\ldots , a_m]}
$$
$$
\cong \bigoplus _{k_1=0}^{n-1}\cdots \bigoplus _{k_m=0}^{n-1} F_{[a_1-k_1,\ldots  , a_m -k_m]} .
$$
\end{lemma}
\begin{proof}
Indeed, we have
$$
\left(\cO (-\sum _{i=1}^m \f{k_i}{n}D_i) \otimes F \right) _{[a_1,\ldots , a_m]}
\cong F_{[a_1-k_1,\ldots  , a_m -k_m]} 
$$
by Lemma \ref{formulae} above. 
\end{proof}

\begin{lemma}
\label{secondcomponents}
$$
\left( \bigoplus _{u_1=0}^{n-1}\cdots \bigoplus _{u_m=0}^{n-1} F_{[u_1,\ldots , u_m]}\otimes 
\cO (-\sum _{i=1}^m \f{u_i}{n}D_i)\right) _{[a_1,\ldots , a_m]}
$$
$$
\cong \bigoplus _{u_1=0}^{n-1}\cdots \bigoplus _{u_m=0}^{n-1} 
F_{[u_1,\ldots , u_m]}\otimes \cO (\sum _{i=1}^m [\frac{a_i - u_i}{n}]D_i).
$$
\end{lemma}
\begin{proof}
We have 
$$
\left( 
\cO (-\sum _{i=1}^m \f{u_i}{n}D_i)\right) _{[a_1,\ldots , a_m]} =
\cO (\sum _{i=1}^m [\f{a_i - u_i}{n}]D_i)
$$
and hence, since $F_{[u_1,\ldots , u_m]}$ is just a vector bundle on $X$,
$$
\left( F_{[u_1,\ldots , u_m]}\otimes 
\cO (-\sum _{i=1}^m \f{u_i}{n}D_i)\right) _{[a_1,\ldots , a_m]} \cong
F_{[u_1,\ldots , u_m]}\otimes
\cO (\sum _{i=1}^m [\frac{a_i - u_i}{n}]D_i).
$$
\end{proof}

We put these two together with the following.

\begin{lemma}
\label{comparecomponents}
$$
\bigoplus _{u_1=0}^{n-1}\cdots \bigoplus _{u_m=0}^{n-1} 
F_{[u_1,\ldots , u_m]}\otimes \cO (\sum _{i=1}^m [\frac{a_i - u_i}{n}]D_i)
$$
$$
\cong \bigoplus _{k_1=0}^{n-1}\cdots \bigoplus _{k_m=0}^{n-1} F_{[a_1-k_1,\ldots  , a_m -k_m]} .
$$
\end{lemma}
\begin{proof}
For given integers $0\leq a_i <n$ and $0\leq u_i < n$, set
$$
k_i := a_i -u_i - n\cdot [\frac{a_i - u_i}{n}],
$$ 
so that 
$$
a_i -k_i = u_i + n\cdot [\frac{a_i - u_i}{n}].
$$
With this definition of $k_i$ we have
$$
F_{[u_1,\ldots , u_m]}\otimes \cO (\sum _{i=1}^m [\frac{a_i - u_i}{n}]D_i)
\cong
F_{[a_1-k_1,\ldots  , a_m -k_m]} ,
$$
due to the periodicity of the parabolic structure.

Note that $0\leq k_i <n$, because  $a_i - u-i < 0$ if and only if the greatest integer piece in the definition of
$k_i$ is equal to $-1$ (otherwise it is $0$). 

For a fixed $(a_1,\ldots , a_m )$, as $(u_1,\ldots , u_m)$ ranges over all possible choices
with $0\leq u_i <n$ the resulting $(k_1, \ldots, k_m)$
also ranges over all possible choices with $0\leq k_i <n$.
Thus we get the isomorphism which is claimed.
\end{proof}

\begin{corollary}
\label{cchiso}
If $F$ is a parabolic bundle on $X$ with respect to $D_1,\ldots , D_m$ then the parabolic bundles
$$
\left( \bigoplus _{k_1=0}^{n-1}\cdots \bigoplus _{k_m=0}^{n-1} \cO (-\sum _{i=1}^m \f{k_i}{n}D_i)\right) \otimes F
$$
and
$$
\bigoplus _{u_1=0}^{n-1}\cdots \bigoplus _{u_m=0}^{n-1} F_{[u_1,\ldots , u_m]}\otimes 
\cO (-\sum _{i=1}^m \f{u_i}{n}D_i)
$$
are componentwise isomorphic, hence componentwise Chow equivalent.
\end{corollary}
\begin{proof}
Putting together Lemmas \ref{firstcomponents}, \ref{secondcomponents} and \ref{comparecomponents} gives,
for any $0\leq a_i <n$ 
$$ 
\left(
\left( \bigoplus _{k_1=0}^{n-1}\cdots \bigoplus _{k_m=0}^{n-1} \cO (-\sum _{i=1}^m \f{k_i}{n}D_i)\right) \otimes F \right)_{[a_1,\ldots , a_m]}
$$
$$
\cong
\left( \bigoplus _{u_1=0}^{n-1}\cdots \bigoplus _{u_m=0}^{n-1} F_{[u_1,\ldots , u_m]}\otimes 
\cO (-\sum _{i=1}^m \f{u_i}{n}D_i)\right) _{[a_1,\ldots , a_m]}.
$$
\end{proof}

We can now calculate with the previous corollary.

\begin{theorem}
\label{answer}
Suppose $F$ is a parabolic bundle on $X$ with respect to $D_1,\ldots , D_m$, with $n$ as denominator. 
Then we have the following
formula for the Chern character of $F$:
$$
{\rm ch}(F) = \frac{\sum _{a_1=0}^{n-1}\cdots \sum _{a_m=0}^{n-1} e^{-\sum _{i=1}^m \f{a_i}{n}D_i}
{\rm ch}(F_{[a_1,\ldots , a_m]})}
{\sum _{a_1=0}^{n-1}\cdots \sum _{a_m=0}^{n-1} e^{-\sum _{i=1}^m \f{a_i}{n}D_i} }.
$$
In other words, the Chern character of $F$ is the weighted average of the Chern characters of 
the component bundles in the range $0\leq a_i < n$, with weights $e^{-\sum _{i=1}^m \f{a_i}{n}D_i}$.
\end{theorem}
\begin{proof}
By Theorem \ref{invariance-th}, two componentwise Chow equivalent parabolic bundles have the same Chern character. From 
the general theory over a DM stack we know that Chern character of parabolic bundles is additive
and multiplicative, and Lemma \ref{formulae} says that it behaves as usual on line bundles. Therefore
the Chern characters of both parabolic bundles appearing in the statement of Corollary \ref{cchiso}
are the same. This gives the formula
$$
(\sum _{a_1=0}^{n-1}\cdots \sum _{a_m=0}^{n-1} e^{-\sum _{i=1}^m \f{a_i}{n}D_i} )\cdot 
{\rm ch}(F) 
= \sum _{a_1=0}^{n-1}\cdots \sum _{a_m=0}^{n-1} e^{-\sum _{i=1}^m \f{a_i}{n}D_i}
{\rm ch}(F_{[a_1,\ldots , a_m]})
$$
The term multiplying ${\rm ch}(F)$ on the left side is an element of the Chow group which has nonzero term 
in degree zero. Therefore, in the rational Chow group it can be inverted and we get the formula stated in the
theorem.
\end{proof}

\begin{remark}
The function
$$
(a_1,\ldots , a_m)\mapsto  e^{-\sum _{i=1}^m \f{a_i}{n}D_i}
{\rm ch}(F_{[a_1,\ldots , a_m]})
$$
is periodic in the variables $a_i$, that is the value for $a_i+n$ is the same as the value for $a_i$.
\end{remark}

\begin{remark}
Also the formula is clearly additive. 
\end{remark}


\section{Examples with parabolic line bundles}


We verify the formula of Theorem \ref{answer} for parabolic line bundles, and then give some examples which are 
direct sums of line bundles which show why it
is necessary to include all of the terms $F_{[a_1,\ldots , a_m]}$ in the formula. 

\subsection{Verification for line bundles}

Suppose $L=\cO(\al.D)$ is a parabolic line bundle on $(X,D)$ where $D$ is an irreducible and 
smooth divisor and $\al=\f{h}{n}\, \in \,\Q$. The formula of Theorem \ref{answer} is obviously invariant if
we tensor the parabolic bundle by a vector bundle on $X$, in particular we can always tensor with an integer power of
$\cO (D)$ so it suffices to check when $0\leq h<n$. 

Notice that the constituent bundles are 
\begin{eqnarray*}
L_{[a_i]} & =  \cO &\m{ if } 0\leq a_i \leq n-h-1  \\
          &\,\,\,\,=  \cO(D) &   \m{ if }n-h\leq  a_i <n. \\
\end{eqnarray*}

We have to check that the Chern character of $L$ is 
$$ {\rm ch} (L)= e^{\al.D}$$

The formula in Theorem \ref{answer} gives

\begin{eqnarray*}
{\rm ch}(L) &=& \f{1 + e^{-\f{1}{n}.D} +...+ e^{-\f{(n-h-1)}{n}.D}+ e^{-\f{(n-h)}{n}.D}.e^D+...
+e^{-\f{(n-1)}{n}.D}.e^D}{1+e^{-\f{1}{n}.D}+...+e^{-\f{(n-1)}{n}.D}}\\
&=& \f{1 + e^{-\f{1}{n}.D} +...+ e^{-\f{(n-h-1)}{n}.D}+ e^{-\f{(n-h)}{n}.D}.e^D+...
+e^{-\f{(n-1)}{n}.D}.e^D} {(\f{1}{e^{\f{h}{n}.D}})(e^{\f{h}{n}.D}+e^{\f{h-1}{n}.D}+...+1+...+e^{-\f{n-1-h}{n}.D})} \\
&=& e^{\f{h}{n}.D}
\end{eqnarray*}

Suppose $D_1,D_2,...,D_m$ are distinct smooth divisors which have normal crossings on $X$. 
Let $L_i=\cO(\al_i.D_i)$ be parabolic line bundles with $\al_i\in \Q$, for $1\leq i\leq m$.
Then the constituent bundles of the tensor product $L:=L_1\otimes L_2\otimes...\otimes L_m$ are
$$
(L_1\otimes L_2\otimes...\otimes L_m)_{[a_1,a_2,...,a_m]}\,=\,(L_1)_{[a_1]}\otimes (L_2)_{[a_2]}\otimes...\otimes (L_m)_{[a_m]}.
$$
and  
$$
{\rm ch}\left((L_1\otimes L_2\otimes...\otimes L_m)_{[a_1,a_2,...,a_m]}\right)
\,=\,{\rm ch}\left((L_1)_{[a_1]}\right).{\rm ch}\left( (L_2)_{[a_2]}\right)...{\rm ch}\left((L_m)_{[a_m]}\right).
$$
The formula in Theorem \ref{answer} is now easily verified for the case when $L$ is a 
parabolic line bundle as above, once it is verified for the parabolic line bundles $L_i$. 
Indeed, the formula in this case is essentially the product of the Chern characters of
$L_i$, for $1\leq i\leq m$.  

\subsection{The case of two divisors and $n=2$}

Suppose we have two divisor components $D_1$ and $D_2$, and suppose the denominator is $n=2$. Then a parabolic bundle may be
written as a $2\times 2$ matrix
$$
F = \left(  
\begin{array}{cc}
F_{[0,0]} & F_{[0,1]} \\
F_{[1,0]} & F_{[1,1]}
\end{array}
\right) .
$$
In particular by Lemma \ref{trivial} we have 
$$
\cO (\frac{D_1}{2}) = 
\left(  
\begin{array}{cc}
\cO  & \cO \\
\cO (D_1) & \cO (D_1)
\end{array}
\right) ,
\;\;\;
\cO (\frac{D_2}{2}) = 
\left(  
\begin{array}{cc}
\cO  & \cO (D_2) \\
\cO  & \cO (D_2)
\end{array}
\right) ,
$$
$$
\cO (D_1) = 
\left(  
\begin{array}{cc}
\cO  (D_1)& \cO (D_1) \\
\cO (D_1) & \cO (D_1)
\end{array}
\right) ,
\;\;\;
\cO (D_1 + \frac{D_2}{2}) = 
\left(  
\begin{array}{cc}
\cO (D_1) & \cO (D_1+D_2) \\
\cO (D_1) & \cO (D_1+D_2)
\end{array}
\right) ,
$$
and
$$
\cO (\frac{D_1}{2} + \frac{D_2}{2}) = 
\left(  
\begin{array}{cc}
\cO  & \cO (D_2) \\
\cO (D_1) & \cO (D_1 + D_2)
\end{array}
\right) .
$$

\subsection{Counterexample for filtrations}

Giving a parabolic bundle by filtrations amounts essentially to considering the bundle $E = F_{[0,0]}$ together with
its subsheaves $F_{[-1,0]}$ and $F_{[0,-1]}$. By the formula \eqref{extensioncond} these subsheaves are determined
by $F_{[1,0]}$ and $F_{[0,1]}$, that is the upper right and lower left places in the matrix. The lower right place 
doesn't intervene in the filtration notations. This lets us construct an example: if 
$$
F := \cO (\frac{D_1}{2}) \oplus \cO (\frac{D_2}{2}), \;\;\; 
G := \cO \oplus \cO (\frac{D_1}{2} + \frac{D_2}{2})
$$
then $F$ and $G$ have the same underlying bundle $E= \cO \oplus \cO$, and the Chern data for their filtrations are the same,
however their Chern characters are different. For example if
 $X=\p ^2$ and $D_1$ and $D_2$ are two distinct lines whose class is denoted $H$
then
\begin{equation}\label{caseC}
\m{ch} (F) = 
\m{ch}\left(\cO _X(\frac{1}{2}D_1 + \frac{1}{2}D_2) \oplus \cO _X\right) 
= 1 + e^{\frac{1}{2}D_1 + \frac{1}{2}D_2} = 2 + H + \frac{H^2}{2}
\end{equation}
and 
\begin{equation}\label{caseD}
\m{ch}(G) =
\m{ch}\left(\cO _X(\frac{1}{2}D_1) \oplus \cO _X(\frac{1}{2}D_2)\right) = e^{\frac{1}{2}D_1} + e^{\frac{1}{2}D_2}
= 2 e^{H/2} = 2 + H + \frac{H^2}{4}.  
\end{equation}


\subsection{Counterexample for MY structure}

Similarly, the MY-parabolic structure consists of $F_{[0,0]}$ and $F_{[1,1]}$, that is the diagonal terms in the matrix,
and the off-diagonal terms don't intervene. A different example serves to show that there is no easy formula for the
Chern character in terms of these pieces only. Put
$$
F:= \cO (\frac{D_1}{2} + \frac{D_2}{2}) \oplus \cO (D_1), \;\;\;
G:= \cO (D_1 + \frac{D_2}{2}) \oplus \cO (\frac{D_1}{2}).
$$
Then
$$
F_{[0,0]} = \cO \oplus \cO (D_1) = G_{[0,0]}
$$
and 
$$
F_{[1,1]} = \cO (D_1 + D_2) \oplus \cO (D_1) = G_{[1,1]}.
$$
On the other hand, again in the example $X=\p ^2$ and $D_1$ and $D_2$ are lines whose class is denoted $H$ we have
$$
\m{ch}(F) = e^{\frac{1}{2}D_1 + \frac{1}{2}D_2} + e^{D_1} = 2e^H = 2 + 2 H + H^2
$$
whereas
$$
\m{ch} (G) = e^{D_1 + \frac{D_2}{2}} + e^{\frac{D_1}{2}} = e^{\frac{3H}{2}} + e^{\frac{H}{2}} = 2 + 2H + \frac{5}{4}H^2.
$$

In both of these examples, of course the structure with filtrations or the MY-parabolic structure permits to obtain back the
full multi-indexed structure and therefore to get the Chern character, however these examples show that the Chern character cannot
be written down easily just in terms of the Chern characters of the component pieces.


\section{A formula involving intersection of filtrations}\label{gysin}


In this section we will give another expression for the parabolic Chern character formula, when the parabolic structure
is viewed as coming from filtrations on the divisor components. 
This formula will involve terms on the multiple intersections of the divisor components, of
intersections of the various filtrations. 
To see how these terms show up in the formula, we first illustrate it by an example below. 


\subsection{Example on surfaces}\label{example}

We consider more closely how the intersection of the filtrations on $D_1$ and
$D_2$ comes into play for determining the Chern character. 
Panov \cite{Panov} and Mochizuki \cite{Mochizuki} considered this situation and obtained formulas for
the second parabolic Chern class involving intersections of the filtrations. 

For this example we keep the hypothesis that $X$ is a surface
and the denominator is $n=2$, also assuming that there are only two divisor components $D_1$ and $D_2$ intersecting at a point $P$.
The typical example is $X = \p ^2$ and the $D_i$ are distinct lines meeting at $P$.

Let $E= F_{[0,0]}$ be a rank two bundle. 
Consider rank one strict subbundles $B_i\subset E|_{D_i}$. Note that
$$
E|_{D_1} = F_{[0,0]} / F_{[-2,0]}, \;\;\; E|_{D_2} = F_{[0,0]} / F_{[0,-2]}.
$$
There is a unique parabolic structure with 
$$
B_1 = F_{[-1,0]} / F_{[-2,0]},
$$
and 
$$
B_2 = F_{[0,-1]} / F_{[0,-2]}.
$$
The quotient $(E|_{D_1})/B_1$ is a line bundle on $D_1$ and similarly for $D_2$, and if the parabolic 
structure corresponds to the $B_i$ as above then
$$
(E|_{D_1})/B_1 = F_{[0,0]}/F_{[-1,0]}
$$
and similarly on $D_2$. 

In particular, $F_{[-1,0]}$ is defined by the exact sequence
$$
0 \rightarrow F_{[-1,0]}\rightarrow E \rightarrow (E|_{D_1})/B_1\rightarrow 0.
$$
Similarly,  $F_{[0,-1]}$ is defined by the exact sequence
$$
0 \rightarrow F_{[0,-1]}\rightarrow E \rightarrow (E|_{D_2})/B_2\rightarrow 0.
$$

Note that the Chern characters of $F_{[-1,0]}$ and $F_{[0,-1]}$ don't depend on the intersection of the 
$B_i$ over $D_1\cap D_2$. On the other hand, $F_{[-1,-1]}$ is a vector bundle, by the locally abelian condition.
Furthermore, as a subsheaf of $E$ it is equal to $F_{[-1,0]}$ along $D_1$  and $F_{[0,-1]}$ along $D_2$.
Thus, in fact $F_{[-1,-1]}$ is the subsheaf of $E$ which is the intersection of these two subsheaves.
To prove this note that the intersection of two reflexive subsheaves of a reflexive sheaf is again reflexive
because it has the Hartogs exension property. In dimension two, reflexive sheaves are vector bundles,
and they are determined by what they are in codimension one.

We have a left exact sequence
$$
0 \rightarrow F_{[-1,-1]} \rightarrow E \rightarrow (E|_{D_1})/B_1\oplus (E|_{D_2})/B_2 .
$$
Here is where the intersection of the filtrations comes in: in our example $D_1\cap D_2$ is a single
point, denote it by $P$. We have one-dimensional subspaces of the two dimensional fiber of $E$ over $P$:
$$
B_{1,P}, B_{2,P}\subset E_P.
$$
There are two cases: either they coincide, or they don't. 

{\em Case A: they coincide}---In this case we can choose a local frame for $E$ in which $B_1$ and $B_2$
are both generated by the first basis vector. We are basically in the direct sum of two rank one bundles, 
one of which containing the two subspaces and the other not. In this case there is an exact sequence
$$
0 \rightarrow F_{[-1,-1]} \rightarrow E \rightarrow (E|_{D_1})/B_1\oplus (E|_{D_2})/B_2 \rightarrow Q \rightarrow 0
$$
where $Q$ is a rank one skyscraper sheaf at $P$.
This is because the fibers of $(E|_{D_1})/B_1$ and $(E|_{D_2})/B_2 $ coincide at $P$,
and $Q$ is by definition this fiber with the map being the difference of the two elements.
Things coming from $E$ go to the same in both fibers so they map to zero in $Q$.

An example of this situation would be the parabolic bundle $\cO _X(\frac{1}{2}D_1 + \frac{1}{2}D_2) \oplus \cO _X$. 

{\em Case B: they differ}---In this case we can choose a local frame for $E$ in which $B_1$ and $B_2$ 
are generated by the two basis vectors respectively. In this case the map in question is surjective
so we get a short exact sequence 
$$
0 \rightarrow F_{[-1,-1]} \rightarrow E \rightarrow (E|_{D_1})/B_1\oplus (E|_{D_2})/B_2 \rightarrow 0.
$$

An example of this situation would be the parabolic bundle $\cO _X(\frac{1}{2}D_1) \oplus \cO _X(\frac{1}{2}D_2)$.

The formula for the Chern character will involve the Chern character of $E$, the Chern characters of the
bundles $B_i$, and a correction term for the intersection.  All other things being equal,
the formulas in the two cases will differ by ${\rm ch}(Q)$ at the place $F_{[1,1]}$ (this is the
same as for $F_{[-1,-1]}$). When the weighted average is taken, this comes in with a coefficient of 
$(\frac{1}{4} +\ldots )$, but the higher order terms multiplied by the codimension $2$ class 
${\rm ch}(Q)$ come out to zero
because we are on a surface. Therefore, the formulae in case A and case B will differ by $\frac{1}{4}
{\rm ch}(Q)$. Fortunately enough this is what actually happens in the examples of the previous section!


\subsection{Changing the indexing}
\label{subsec-changing}

When describing a parabolic bundle by filtrations, we most naturally get to the bundles $F_{[a_1,\ldots , a_m]}$
with $-n \leq a_i \leq 0$. On the other hand, the weighted average in Theorem \ref{answer} is over 
$a_i$ in the positive interval $[0,n-1]$. It is convenient to have a formula which brings into play the
bundles in a general product of intervals. The need for such was seen in the example of the previous subsection. 

We  have the following result which meets up with Mochizuki's notation and discussion in \cite[\S 3.1]{Mochizuki}.

\begin{proposition}
\label{change}
Let ${\bf b}= (b_1,\ldots , b_m)$ be any multi-index of integers. Then ${\rm ch}(F)$ is obtained by taking the
weighted average of the ${\rm ch}(F_{[a_1,\ldots , a_m]})$ with weights $ e^{-\sum _{i=1}^m \f{a_i}{n}D_i}$,
over the product of intervals $b_i \leq a_i < b_i+n$, and then multiplying by $e^{-\sum _{i=1}^m \f{b_i}{n}D_i}$
(that is the weight for the smallest multi-index in the range). This formula may also be written as:
\begin{equation}
\label{extraformula}
{\rm ch}(F) = \frac{\sum _{a_1=b_1}^{b_1+n-1}\cdots \sum _{a_m=b_m}^{b_m+n-1} e^{-\sum _{i=1}^m \f{a_i}{n}D_i}
{\rm ch}(F_{[a_1,\ldots , a_m]})}
{\sum _{a_1=0}^{n-1}\cdots \sum _{a_m=0}^{n-1} e^{-\sum _{i=1}^m \f{a_i}{n}D_i} }.
\end{equation}
\end{proposition}
\begin{proof}
If $a_i$ and $a'_i$ differ by integer multiples of $n$ then by using condition \eqref{extensioncond} of \S \ref{se.-parfilt}, we have
$$
e^{-\sum _{i=1}^m \f{a'_i}{n}D_i}
{\rm ch}(F_{[a'_1,\ldots , a'_m]}) =
e^{-\sum _{i=1}^m \f{a_i}{n}D_i}
{\rm ch}(F_{[a_1,\ldots , a_m]}).
$$
Thus, the numerator in the formula \eqref{extraformula} is equal to the numerator of the formula
in Theorem \ref{answer}. The denominators are the same. On the other hand, if we form the weighted average as described
in the first sentence of the proposition, then the numerator will be the same as in \eqref{extraformula}.
The denominator of the weighted average is
$$
\sum _{a_1=b_1}^{b_1+n-1}\cdots \sum _{a_m=b_m}^{b_m+n-1} e^{-\sum _{i=1}^m \f{a_i}{n}D_i} = 
e^{-\sum _{i=1}^m \f{b_i}{n}D_i}\sum _{a_1=0}^{n-1}\cdots \sum _{a_m=0}^{n-1} e^{-\sum _{i=1}^m \f{a_i}{n}D_i}.
$$
Hence, when we multiply the weighted average by $e^{-\sum _{i=1}^m \f{b_i}{n}D_i}$ we get \eqref{extraformula}. 
\end{proof}

\begin{remark}
If we replace the denominator $n$ by a new one $np$ then 
the formulae of Theorem \ref{answer} or the previous proposition, give the same answers.
\end{remark}
Indeed, the parabolic structure $\widetilde{F}$ for denominator $np$ contains the same 
sheaves, but each one is copied $p^m$ times:
$$
\widetilde{F}_{[pa_1 + q_1 , \ldots , pa_m + q_m]} = F_{[a_1,\ldots , a_m]}
$$
for $0 \leq q_i \leq p-1$. Therefore,  
both the numerator and the denominator in our formulae are multiplied by 
$$
\sum _{q_1=0}^{p-1}\cdots \sum _{q_m=0}^{p-1} e^{-\sum _{i=1}^m \f{q_i}{np}D_i},
$$
and the quotient stays the same.


\subsection{A general formula involving intersection of filtrations}
\label{subsec-gysin}

We can generalize the example of surfaces  in \S \ref{example}, to get a formula which generalizes
the codimension 2 formulae of Panov \cite{Panov} and Mochizuki \cite{Mochizuki}. 

In this section we suppose we are working with the notation of a locally abelian parabolic structure $F$ given by filtrations,
on a vector bundle $E:= F_{[0,\ldots , 0]}$ with filtrations  
$$
E|_{D_i}=F^i_{0}\supset F^i_{-1}\supset...\supset F^i_{-n}=0.
$$
Then for $-n \leq a_i \leq 0$ define the quotient sheaves supported on $D_i$
$$
Q^i_{[a_i]} := \f{E|_{D_i}}{F^i_{a_i}} 
$$
and the parabolic structure $F_{\cdot}$ is given by 
\begin{equation}
\label{Fformula}
F_{[a_1,\ldots , a_m]} = \ker \left( E \rightarrow \oplus _{i=0}^m Q^i_{[a_i]} \right) .
\end{equation}

More generally define a family of multi-indexed quotient sheaves by
\begin{eqnarray*}
Q^i_{[a_i]}& := & \f{E|_{D_i}}{F^i_{a_i}} \m{ on } D_i \\
Q^{i,j}_{[a_i,a_j]}& := & 
\f{E|_{D_i\cap D_j}}{F^i_{a_i} + F^j_{a_j}}
\m{ on }D_i\cap D_j \\
& \vdots & \\
Q_{[a_1,a_2,...,a_m]} & := & 
\f{E|_{D_1\cap \ldots \cap D_m}}{F^1_{a_1} + \ldots + F^m_{a_m}}
\m{ on } D_1\cap D_2\cap \ldots \cap D_m.
\end{eqnarray*}
In these notations we have $-n \leq a_i \leq 0$. 

If we consider quotient sheaves as corresponding to linear subspaces of the Grothendieck projective bundle
associated to $E$, then the multiple quotients above are multiple intersections of the $Q^i_{[a_i]}$. 
The formula \eqref{Fformula} extends to a Koszul-style resolution of the component sheaves of the parabolic structure.

\begin{lemma}
\label{resolution}
Suppose that the filtrations give a locally abelian parabolic structure,
in particular they satisfy the conditions of Lemma \ref{locabconseq}. Then
for any $-n \leq a_i \leq 0$
the following sequence is well defined and exact: 
$$
0\rar F_{[a_1,a_2,...,a_m]}\rar E \rar 
\bigoplus_{i=1}^m Q^i_{[a_i]}  
\rar \bigoplus_{i<j}Q^{i,j}_{[a_i,a_j]} \rar ...\rar
Q_{[a_1,a_2,...,a_m]}\rar 0.
$$
\end{lemma}
\begin{proof}
The maps in the exact sequence are obtained from the quotient structures of the terms
with alternating signs like in the \^Cech complex. 
We just have to prove exactness. This is a local question. 
By the locally abelian condition, we may assume that
$E$ with its filtrations is a direct sum of rank one pieces. The formation of the sequence, and
its exactness, are compatible with direct sums. Therefore we may assume that $E$ has rank one,
and in fact $E\cong\cO _X$. 

In the case where $E$ is the rank one trivial bundle,
the filtration steps are either $0$ or all of $\cO_{D_i}$. In particular, there is 
$-n < b_i\leq 0$ such that $F^i_{j} = \cO_{D_i}$ for $j\geq b_i$ and $F^i_{j} = 0$ for $j < b_i$. 
Then
$$
Q^{i_1,\ldots ,i_k}_{[a_{i_1},\ldots ,a_{i_k}]} = \cO_{D_{i_1}\cap \cdots \cap D_{i_k}}
$$
if $a_{i_j}< b_{i_j}$ for all $j= 1,\ldots , k$, and the quotient is zero otherwise. 

The sequence is defined for each multiindex $a_1,\ldots , a_m$. Up to reordering the coordinates
which doesn't affect the proof, we may assume that there is $p\in [0,m]$ such that $a_i < b_i$ for $i\leq p$,
but $a_{i}\geq b_i$ for $i>p$.  In this case, the quotient is nonzero only when $i_1,\ldots , i_k \leq p$. 
Furthermore, 
$$
F_{[a_1,\ldots ,a_m]} = \cO (-D_1-\cdots -D_p).
$$
In local coordinates, the divisors $D_1,\ldots , D_p$ are coordinate divisors. Everything is constant in
the other coordinate directions which we may ignore.
The complex in question becomes 
$$
\cO (-D_1-\cdots -D_p) \rar \cO \rar \oplus_{ 1\leq i \leq p} \cO _{D_i} \rar
\oplus _{1\leq i < j \leq p} \cO _{D_i\cap D_j} \rar \ldots \rar \cO _{D_1\cap \cdots \cap D_p} . 
$$
Etale locally, this is exactly the same as the exterior tensor product of $p$ copies of the
resolution of $\cO _{\A^1}(-D)$ on the affine line $\A^1$ with divisor $D$ corresponding to the origin,
$$
\cO _{\A ^1}(-D)\lrar \cO _{\A^1} \lrar \cO _D \lrar 0.
$$
In particular, the exterior tensor product complex is exact except at the beginning, where it resolves
$\cO (-D_1-\cdots -D_p)$ as required.
\end{proof}

Using the resolution of Lemma \ref{resolution} we can compute the Chern character of $F_{[a_1,a_2,...,a_m]}$ in
terms of the Chern character of sheaves supported on intersection of
the divisors $D_{i_1}\cap...\cap D_{i_r}$. This gives us
$$
{\rm ch}(F_{[a_1,a_2,...,a_m]}) \,=\, {\rm ch}(E) + \sum_{k=1}^m (-1)^k \sum_{i_1<i_2<...<i_k}{\rm ch}(Q_{[a_{i_1},\ldots , a_{i_k}]}^{i_1,...,i_k}).
$$

Substituting this formula for ${\rm ch}(F_{[a_1,a_2,...,a_m]})$ into
Theorem \ref{answer}, or rather into \eqref{extraformula} of Proposition \ref{change} with $b_i = -n$, 
we obtain the following formula for the
associated parabolic bundle. Note that the limits of the sums are different
in the numerator and denominator, as in \eqref{extraformula}. Also the term ${\rm ch}(E)$ occurs with a different factor
in the numerator and denominator; the ratio of these factors is $e^{\sum _{i=1}^mD_i}=e^D$. 

\begin{corollary}\label{co.-gysinformula}
If $F$ is a locally abelian parabolic bundle then
$$
{\rm ch}(F) = e^D{\rm ch}(E) + 
$$
$$
\frac{\sum _{a_1=-n}^{-1}\cdots \sum _{a_m=-n}^{-1} e^{-\sum _{i=1}^m \f{a_i}{n}D_i}
\sum_{k=1}^m (-1)^k \sum_{i_1<i_2<...<i_k}{\rm ch}(Q_{[a_{i_1},\ldots , a_{i_k}]}^{i_1,...,i_k})}
{\sum _{a_1=0}^{n-1}\cdots \sum _{a_m=0}^{n-1} e^{-\sum _{i=1}^m \f{a_i}{n}D_i} }.
$$
\end{corollary}

In fact, we can also write the formula in terms of an associated graded. For this, fix $1\leq i_1 < \cdots < i_k \leq m$
and analyze the quotient $Q ^{i_1,\ldots ,i_k}_{[a_{i_1}-1,\ldots ,a_{i_k}-1]}$ along the multiple intersection
$D_{i_1\cdots i_k}$. There, the bundle $E|_{D_{i_1\cdots i_k}}$ has  $k$ filtrations $F^{i_j}_{a_{i_j}}|_{D_{i_1\cdots i_k}}$
indexed by $-n \leq a_{i_j}\leq 0$,
leading to a multiple-associated-graded 
defined as follows. 
For $-n \leq a_{i_j}\leq 0$ put
$$
F ^{i_1,\ldots ,i_k}_{[a_{i_1},\ldots ,a_{i_k}]} := \bigcap _{j=1}^k F^{i_j}_{a_{i_j}}|_{D_{i_1\cdots i_k}} .
$$
Then define
\begin{equation}
\label{multigraded}
Gr ^{i_1,\ldots ,i_k}_{[a_{i_1},\ldots ,a_{i_k}]} := 
\frac{F ^{i_1,\ldots ,i_k}_{[a_{i_1},\ldots ,a_{i_k}]}}{\sum _{q=1}^k F ^{i_1,\ldots ,i_k}_{[a_{i_1},\ldots ,a_{i_q}-1, \ldots a_{i_k}]}}
\end{equation}
where the indices in the denominator are almost all $a_{i_j}$ but one $a_{i_q}-1$. 
A good way to picture this when $k=2$ is to draw a square divided into a grid whose sides are the intervals $[-n,0]$. The filtrations correspond to
horizontal and vertical half-planes intersected with the square. Pieces of the
associated-graded are indexed by grid squares, indexed by their upper right points. Thus the pieces are defined for 
$1-n \leq a_{i_j} \leq 0$.

If the parabolic structure is locally abelian then the filtrations admit a common splitting and we have
$$
Gr ^{i_1,\ldots ,i_k}_{[a_{i_1},\ldots ,a_{i_k}]} = 
Gr ^{F^{i_1}}_{a_{i_1}}Gr ^{F^{i_2}}_{a_{i_2}} \cdots Gr ^{F^{i_k}}_{a_{i_k}}(E|_{D_{i_1\cdots i_k}}),
$$
or more generally the same thing in any order. Without the common splitting hypothesis, the multi-graded defined previously would not
even have dimensions which add up. 

The multi-quotient has an induced multiple filtration whose associated-graded is a sum of pieces of the multi-graded defined above. 
In the $k=2$ picture, the multi-quotient corresponds to a rectangle in the upper right corner of the square. 
For example, we have
$$
Gr ^{i_1,\ldots ,i_k}_{[a_{i_1},\ldots ,a_{i_k}]} \cong 
\ker \left(  
Q ^{i_1,\ldots ,i_k}_{[a_{i_1}-1,\ldots ,a_{i_k}-1]} 
\rightarrow \bigoplus _{j=1}^k
Q ^{i_1,\ldots ,i_k}_{[a_{i_1}-1,\ldots , a_{i_j},\ldots  ,a_{i_k}-1]} 
\right)
$$
where in the direct sum, the indices are all $a_{i_l}-1$ except for one which is $a_{i_j}$.

Thus in the Grothendieck group of sheaves on $D_{i_1}\cap \cdots \cap D_{i_k}$, we have an equivalence
$$
Q ^{i_1,\ldots ,i_k}_{[a_{i_1},\ldots ,a_{i_k}]}\sim 
\bigoplus _{\underline{c},\,  a_{i_j}< c_{i_j} \leq 0}
Gr ^{i_1,\ldots ,i_k}_{[c_{i_1},\ldots ,c_{i_k}]} .
$$
This gives us the following formula, based on Corollary 
\ref{co.-gysinformula} which in turn comes from \eqref{extraformula} of
Proposition \ref{change} (thus as before the limits of the sum in the numerator and denominator are different). 

\begin{corollary}\label{co.-gysinformulagraded}
Suppose $F$ is a locally abelian parabolic structure. Define the multi-associated-graded by \eqref{multigraded} above. 
Then we have the formula 
{\small
$$
{\rm ch}(F)= 
e^D {\rm ch}(E) \,+\, 
$$
$$
\frac{\sum _{-n \leq a_1,\ldots , a_m <0}e^{-\sum _{i=1}^m \f{a_i}{n}D_i}
\sum_{k=1}^m (-1)^k \sum_{i_1<i_2<...<i_k} \sum _{a_{i_j}< c_{i_j} \leq 0}
{\rm ch} (Gr ^{i_1,\ldots ,i_k}_{[c_{i_1},\ldots ,c_{i_k}]})}
{\sum _{a_1=0}^{n-1}\cdots \sum _{a_m=0}^{n-1} e^{-\sum _{i=1}^m \f{a_i}{n}D_i} }.
$$}
\end{corollary}


\subsection{The case of a single smooth divisor}
\label{subsec.-singlesmooth}

In the case when there is only one smooth divisor component $D$ this formula becomes
\begin{equation}
\label{onecomponent}
{\rm ch}(F) =  
e^D{\rm ch}(E) -
\frac{\sum _{-n < c \leq 0} \left( \sum _{-n \leq a <c}e^{-\frac{a}{n}D} \right) {\rm ch}(Gr_{[c]})}{\sum _{0\leq a <n} e^{-\frac{a}{n}D}} . 
\end{equation}
This can be simplified using the identity $(1+x+ \ldots + x^{n-1}) = (1-x)^{-1}(1-x^n)$ applied to
$x=e^{-\frac{1}{n}D}$, which gives
$$
{\rm ch}(F) =  
e^D {\rm ch}(E) -
\sum _{-n< c \leq 0}  \frac{e^D -e ^{-\frac{c}{n}D}}{1-e^{-D}}{\rm ch}(Gr_{[c]}). 
$$
We can again rewrite this in terms of the rational indexing in the interval $(-1,0]$, denoting by $Gr_{\alpha}$ the
graded $Gr_{[n\alpha ]}$. The formula becomes
\begin{equation}
\label{onecomponent2}
{\rm ch}(F) =  
e^D {\rm ch}(E) -
\sum _{-1< \alpha \leq 0}  \frac{e^D -e ^{-\alpha D}}{1-e^{-D}}{\rm ch}(Gr_{\alpha}). 
\end{equation}
It should be interpreted formally, in the sense that the fractions are
first reduced to power series, then applied to $D\in CH^1(X)_{\Q}$, resulting in polynomials 
because of the nilpotence of the product structure
on $CH^{>0}(X)_{\Q}$. 

One checks that it gives the right formula for a line bundle $F=\cO (\frac{b}{n}D)$. 
We leave it to the reader to make the analogous transformations of the formula in the case of several divisors, possibly
meeting only pairwise as a start, 
and to compare the result with the codimension $2$ formulae of Panov \cite{Panov}
and Mochizuki \cite{Mochizuki}.


\section{Parabolic bundles with real weights }


In this section we consider parabolic bundles with real weights and
define their Chern character and pullback bundles.

Let $X$ be a smooth variety and $D$ be a normal crossing divisor on $X$. Write
$D=D_1+\ldots +D_m$ where $D_i$ are the irreducible smooth components and meeting transversally.

A {\em parabolic bundle} on $(X,D)$ is a collection of vector bundles $F_{\alpha}$  indexed by
multi-indices $\alpha =(\alpha _1,\ldots , \alpha _k)$ with $\alpha _i\in \R$, satisfying the same conditions as
recalled in \S 2. The structure is determined by the sheaves $F_{\alpha}$ for
a finite collection of indices $\alpha$ with $0\leq \alpha _i < 1$, the {\em weights}.

\begin{remark}\label{re.-realwt}
A parabolic bundle with rational weights and denominator $n$
can be considered as a
parabolic bundle with real weights by setting
$$
F_{(t_1,t_2,...,t_m)}:= F_{[[nt_1],[nt_2],...,[nt_m]]} = F_{(\frac{[nt_1]}{n},\frac{[nt_2]}{n},...,\frac{[nt_m]}{n})}
$$
where $[nt_i]$ is the greatest integer less than or equal to $nt_i$, for
any $t_i \in [0,1)\subset \R$.
\end{remark}

We say that $F$ is \textit{locally abelian} if in a Zariski neighbourhood of
any point $x\in X$, $F$ is isomorphic to a direct sum of parabolic
line bundles with real coefficients.


\subsection{Perturbation of parabolic bundles  with real weights}

The following construction is a simplified version of the one considered by Mochizuki \cite[\S3.3]{Mochizuki}, 
and which suffices for our purpose. 
Variations of parabolic weights were considered earlier in \cite{MehtaSeshadri}, \cite{Boden2}, \cite{Thaddeus}.

Suppose $F$ is a parabolic bundle wih real weights on a smooth variety
$(X,D)$.
Consider the real weights
$$
\{\al=(\al_1,\al_2,...,\al_m) :\,0\leq \al_i\leq 1\} .
$$
By definition
$$
F_{\al}|_{D_i}=\f{F_\al}{F_{\al-\delta_i}}
$$
and denote the image
$$
\ov{F_{\al ; D_i,\gamma _i}} := \rm{Im}\left( F_{(\alpha _i,\ldots , \gamma _i, \ldots , \alpha _m )}\lrar F_{\al}|_{D_i} \right)
$$
whenever $\al _i -1 <\gamma _i\leq \al _i$.

Note that if $\gamma$ is a multiindex with 
$\alpha _i -1 < \gamma _i \leq \alpha _i$ then we have an exact sequence
$$
0\rightarrow F_{\gamma } \rightarrow F_{\alpha} \rightarrow \bigoplus _i \frac{F_{\al}|_{D_i}}{\ov{F_{\al ; D_i,\gamma _i }}} .
$$

Consider the graded sheaves 
$$
{\bf{gr}}^i_{\al ;  \gamma_i}F:=
\f{F_{\al}|_{D_i}}{\ov{F_{\al ; D_i,\gamma_i}}}.
$$ 
By the semicontinuity condition there are finitely many indices
and $\gamma_i$ such that 
the graded sheaves ${\bf{gr}}^i_{\al_i-\gamma_i}F$ are non-zero.

Let 
$$ 
r_{\al_i}=\rm{min}\{|\al_i-\gamma_i|\,:\, {\bf{gr}}^i_{\al_i/ \gamma_i}F\neq 0\}
$$

Choose $\ep_{\al_i}$ such that $\ep_{\al_i}<r_{\al_i}$ and $\al_i+\ep_i$ is a rational
number, for each $i$.

The following  construction was used by Mochizuki in \cite[\S 3.4]{Mochizuki}. 

\begin{definition}
A parabolic bundle $F^\ep$ with rational weights
$a_i=\al_i+\ep_{\al_i}$ is defined by setting :
$$
F^\ep_{[a_1,a_2,...,a_m]}:=
F_{\al_1+\ep_{\al_1},..., \al_m+\ep_{\al_m}}.
$$
We call $F^\ep$ as an $\ep$--perturbation of $F$ on $X$.
\end{definition}

For any rational weights $t=[t_1,...,t_m]$, we have the
inclusion of sheaves
$$
F_t\hookrightarrow F^\ep_t 
$$
In other words, we can write
$$
F\hookrightarrow F^\ep.
$$
Write $\ep=\{\ep_{\al_i}\}$, where $\al_i$ runs over the finite set of
real weights which determine $F$.

Suppose $\{ F^i\} _{i\in I}$ is a projective system of parabolic bundles indexed by an ordered set $I$
with inclusions $F^i \hookrightarrow F^j$ for $i\leq j$. Define the {\em intersection}
by the formula
$$
\left( \bigcap_{i\in I} F^i \right) _{\alpha} := \bigcap_{i\in I} F^i _{\alpha}.
$$
This defines a parabolic sheaf. We say that the collection $\{ F^i\} _{i\in I}$ is {\em simultaneously locally
abelian} if there is an etale covering of $X$ such that on the pullback to this etale covering, each of the
$F^i$ admits a direct sum decomposition as a sum of parabolic line bundles, and the inclusion maps are compatible
with these direct sum decompositions. Inclusions of parabolic line bundles are just inequalities of real divisors,
and the intersection of a family of parabolic line bundles just corresponds to taking the $inf$ of the family of real coefficients.
Thus we have the following useful fact.

\begin{lemma}
If $\{ F^i\} _{i\in I}$ is a simultaneously locally abelian projective system of inclusions of parabolic bundles, then
the intersection $\bigcap_{i\in I} F^i$ is a locally abelian parabolic bundle. 
\end{lemma}

\begin{lemma}\label{le.-locab}
Suppose $F$ is a locally abelian parabolic bundle with real weights
$\al=(\al_1,...,\al_m)$ on $(X,D)$. Then any $\ep$--perturbation
$F^\ep$ of $F$ is also locally abelian with the same decomposition. Thus 
the family of $F^{\ep}$ is a simultaneously locally abelian projective system
of inclusions. Taking the intersection we have
$$
F=\bigcap_{\ep\rar 0} F^\ep.
$$
\end{lemma}
\begin{proof}
Since this is a local question, we assume that 
$$
F=\oplus_j \cO(\sum \gamma_i^j.D_i)^{n_j}
$$
for some $\gamma_j^i\in \R$.
Any $\ep$-perturbation of $F$ is
$$
F^\ep= \oplus_j \cO(\sum a^i_j.D_i)^{n_j}
$$
where $a^i_j=\gamma^i_j +\ep^i_j$ are rational numbers and $\ep^i_j$
are small.
Hence $F^\ep$ is locally abelian.
\end{proof}


\subsection{Pullback of parabolic bundles with real weights}

Consider a morphism
$$f:(Y,D')\lrar (X,D)$$
such that $f^{-1}(D)\subset D'$. Here $X,\,Y$ are smooth varieties and
$D,D'$ are normal crossing divisors on $X$ and $Y$ respectively.

In \cite[Lemma 2.6]{Iy-Si}, the pullback of a locally abelian
parabolic bundle with rational weights was defined, using its correspondence with usual
vector bundles on a DM--stack.  
Our aim here is to define the pullback $f^*F$ on $(Y,D')$ of a locally abelian
parabolic bundle $F$ with real weights on $(X,D)$.

\begin{lemma}
\label{le.-pullbackdef}
Suppose $F$ is a locally abelian parabolic bundle with real weights on
$(X,D)$. For any morphism $f:(Y,D')\lrar (X,D)$ such that
$f^{-1}D\subset D'$, we can define the
pullback $f^*F$ on $(Y,D)$ as a locally abelian parabolic bundle with
real weights.
\end{lemma}
\begin{proof}
By Lemma \ref{le.-locab}, we can write
$$
F= \bigcap_{\ep\rar 0}F^\ep.
$$

By \cite[Lemma 2.6]{Iy-Si}, $f^*F^\ep$ is a locally abelian parabolic
bundle with rational weights.

Locally, by Lemma \ref{le.-locab},  each $F^\ep$ is locally abelian,
and the decompositions are compatible for different $\ep$. Thus we can write locally
$$
f^*F^\ep\,=\,\oplus_j \cO(\sum a^i_j(\ep).D'_i)^{n_j}
$$
where $a^i_j(\ep)$ are rational numbers depending on $\ep$.
In other words, the pullbacks form a simultaneously locally abelian projective system.
By Lemma \ref{le.-locab}, we can define the pullback of $F$ as the intersection
$$
f^*F:= \bigcap_{\ep\rar 0}f^*F^\ep ,
$$
and it is a locally abelian parabolic bundle. In fact, locally 
let $\al^i_j=\rm{lim}_{\ep\rar 0}\,a^i_j(\ep)$ (which converges and is
a real number), then 
$$
f^*F\,=\,\oplus_j \cO(\sum a^i_j.D'_i)^{n_j}.
$$
\end{proof}

\subsection{Tensor products of parabolic bundles with real weights}

Suppose $F$ and $G$ are two locally abelian parabolic bundles with real weights. We would like to define their
tensor product. Recall that by 
\cite[Lemma 2.3]{Iy-Si}, the tensor product of locally abelian
parabolic bundles with rational weights can be defined using the correspondence with usual
vector bundles on a DM--stack.

\begin{lemma}
Suppose $F$ and $G$ are locally abelian parabolic bundles with real weights on
$(X,D)$. Then we can define $F\otimes G$ as a locally abelian parabolic bundle with
real weights.
\end{lemma}
\begin{proof}
By Lemma \ref{le.-locab}, we can write
$$
F= \bigcap_{\ep\rar 0}F^\ep, \;\;\; G= \bigcap_{\ep\rar 0}G^\ep
$$
The families $\{ F^{\ep}\}_{\ep \rightarrow 0} $ and $\{ G^{\ep '}\} _{\ep '\rightarrow 0} $ are simultaneously locally abelian, and we can take a
common refinement of the two coverings so that they are locally abelian with respect to the same covering.
Then the family of tensor products $\{ F^{\ep} \otimes G^{\ep '}\} _{\ep , \ep ' \rightarrow 0}$ is again simultaneously locally abelian 
with respect to the same decomposition and we can define
$$
F\otimes G := \bigcap _{\ep , \ep ' \rightarrow 0} F^{\ep} \otimes G^{\ep '}.
$$
\end{proof}

One can also consider duals and internal $\underline{Hom}$.

\subsection{Description by filtrations on a linear constructible decomposition of the space of weights}

For both of the operations defined above, the description in terms of filtrations can jump when the parabolic weights cross ``walls''.
Fix a vector bundle $E$ and filtrations of $E_{D_i}$. These filtrations determine an open subset of possible assignments of weights $\al ^j_i$
to the
filtrations $F^j_i$ with $\al ^{j-1}_i < \al _i^j$. This defines an open subset $W(E,\{ F^j_i\} )\subset \R ^N$. 
Note that the locally abelian condition doesn't depend on the choice of weights but is just a statement about the
filtrations. However, when we apply the pullback operation for a map $(Y,D')\rightarrow (X,D)$ the filtrations on the pullback bundle
might depend on the choice of weights $\underline{\al} \in W(E,\{ F^j_i\} )$. 

A subset of $\R ^N$ is {\em linear-constructible} if it is defined by a finite number of linear equalities and inequalities. 
It is {\em $\Q$-linear-constructible} if the equalities and inequalities have coefficients in $\Q$.

The filtrations for the pullback parabolic bundle are fixed over a $\Q$-linear constructible stratification of the space of weights.
This phenomenon is somewhat similar to what was observed by Budur in \cite{Budur}. 

\begin{proposition}
\label{Qlinearpullback}
Suppose $f:(Y,D')\rightarrow (X,D)$ is a morphism of smooth varieties with normal crossings divisors in good position.
Suppose $(E,\{ F^j_i\} )$ is a locally abelian datum of filtrations for a parabolic structure on $(X,D)$.
There is a stratification of $W(E,\{ F^j_i\} )$ into a finite
disjoint union of  $\Q$-linear constructible sets $W(p)$
such that over each stratum, there is a fixed collection of filtrations $\widetilde{F}^j_i(p)$ for the pullback bundle
$\widetilde{E}:= f^{\ast}E$ and a $\Q$-linear function of weights $f^{\ast}(p):W(p) \rightarrow W(\widetilde{E},\{ \widetilde{F}^j_i(p)\} )$
such that for $\alpha \in W(p)$ the pullback of the parabolic bundle $(E, \{ F^j_i\} , \alpha )$ is equal to 
$(f^{\ast} E, \{ \widetilde{F}^j_i(p)\} , f^{\ast}(p)(\alpha ))$. 
\end{proposition}

We leave the proof to the reader.

A similar statement holds for tensor product, which is again left to the reader. 

\subsection{Chern character of parabolic bundles with real weights}

Suppose $K\subset \R$ is a subfield, and suppose $V$ is a $K$-vector space. If $f\in V\otimes K[x]$ then
we can define in a formal way $\int _0^1f \in V$. 
The same is true if $f$ is a formal piecewise polynomial function whose intervals of different definitions are
defined over $K$. 
A similar remark holds for multiple integrals---in the case we shall consider the domains of piecewise definition will be 
products of intervals defined over $K$ but this could also extend to $K$-linear constructible regions. 

Using this meaning, the formula of Theorem \ref{answer} 
may be rewritten replacing sums by integrals:
\begin{equation}
\label{integralformula}
{\rm ch} (F) = \frac{\int _{\alpha _1=0}^{1}\cdots \int _{\alpha_m=0}^{1} e^{-\sum _{i=1}^m \alpha _iD_i}
{\rm ch}(F_{\alpha})}{\int _{\alpha _1=0}^{1}\cdots \int _{\alpha_m=0}^{1} e^{-\sum _{i=1}^m \alpha _iD_i }}.
\end{equation}
In this formula note that the exponentials of real combinations of divisors are interpreted as formal polynomials.
The power series for the exponential terminates because the product structure of $CH^{>0}(X)$ is nilpotent.

If $F$ is a parabolic bundle with rational weights, then this still takes values in $CH^{\cdot}(X)_{\Q}$. 

If $F$ is a parabolic bundle with real weights, then the formula \eqref{integralformula} may be taken as the
{\em definition} of ${\rm ch}(F)\in CH^{\cdot}(X)_{\R}:= CH^{\cdot}(X) \otimes _{\Z} \R$. No topology 
or metric structure is needed on $CH^{\cdot}(X)_{\R}$ because the integrals involved are piecewise polynomials.

\begin{theorem}
\label{realfunctorial}
The Chern character of locally abelian parabolic bundles with real weights, is additive for
exact sequences, multiplicative for tensor products, and functorial for pullbacks along
good morphisms of varieties with normal crossings divisors. 
\end{theorem}
\begin{proof}
Additivity for exact sequences follows from the shape of the formula. 
Suppose $f: (Y,D')\rightarrow (X,D)$ is a good morphism of varieties with normal crossings divisors.
Fix a bundle and collection of filtrations $(E,\{ F^j_i\} )$ on $(X,D)$. The Chern character may then be 
viewed as a
function 
$$
{\rm ch} : W (E,\{ F^j_i\} )\rightarrow CH^{\cdot}(X)_{\R}.
$$
This function is obtained as a polynomial with coefficients which are rational linear combinations of
the various Chern classes of the intersections of the filtrations, see \S \ref{subsec-gysin}. 
The same may be said of the Chern character of parabolic bundles over $(Y,D')$ once filtrations are fixed. 
Use Proposition \ref{Qlinearpullback} to decompose the space $ W (E,\{ F^j_i\} )$ into a finite union of $\Q$-linear
constructible subsets on which the filtrations of the pullback parabolic structure will be invariant.
Over these subsets the Chern character of the pullback parabolic structures are again 
polynomials with coefficients in $CH^{\cdot}(X)_{\Q}$. On the other hand, by 
\cite[Lemma 2.8]{Iy-Si}, whenever the weights are rational we have that the Chern character of the pullback
is the pullback of the Chern character. We therefore have two polynomials with $CH^{\cdot}(X)_{\Q}$
coefficients which agree on the rational points of a certain $\Q$-linear constructible set.
It follows that the polynomial functions into $CH^{\cdot}(X)_{\R}$ agree on the real points of the
$\Q$-linear constructible set. This proves compatibility of the Chern character for pullbacks of 
real parabolic bundles.

The proof for tensor products is similar, using the analogue of Proposition \ref{Qlinearpullback}. 
\end{proof}


\section{Variants}


In this section we consider a variant of the notion of parabolic structures for the case
of a divisor with multiple points, and also a variant of the construction of 
parabolic bundle associated to a logarithmic connection, concerning the case of unipotent monodromy
at infinity. In both cases, we will restrict to the case when $X$ is a smooth projective surface. 


\subsection{Parabolic structures at multiple points}

Let $X$ be a nonsingular projective surface. Let $D\subset X$ be a divisor such that
$D=\cup_{i=1}^m D_i$ and $D_i$ are smooth and irreducible curves.  Let $P=\{ P_1,\ldots , P_r\}$
be a set of points. 
Assume that the points $P_j$ are crossing points of $D_i$, and that they are general multiple points,
that is through a crossing point $P_j$ we have
divisors $D_1,\ldots , D_k$ which are pairwise transverse. Assume that $D$ has normal crossings outside
of the set of points $P$.

Let $\pi:X'\lrar X$ be the blow--up of $X$ at $P$ and $E$ be the exceptional divisor on
$X'$; note that $E$ is a sum of disjoint exceptional components
$E_j$ over the points $P_j$ respectively.  
The pullback divisor $D'=\sum_{i=1}^m D_i' + E$ is a normal crossing divisor, where $D_i'$ is the 
strict transform of $D_i$, for $1\leq i\leq m$. 

We will define a notion of {\em exceptionally constant parabolic structure on $(X,D,P)$}. 
The term ``exceptionally constant'' means that the parabolic structure pulls back to one which is
constant along the exceptional divisors. Following notation of Mochizuki \cite{Mochizuki} we fix
an origin for the filtrations which is a multi-index ${\bf c}$. This may be important in the present
case since the structures might differ for different values of ${\bf c}$. 

\begin{definition}
Fix a positive integer $n$ for the denominator, and an uplet of
integers ${\bf c} = (c_{D,1},\ldots , c_{D,m}, c_{P,1},\ldots , c_{P,r})$. 
An {\em exceptionally constant parabolic structure on $(X,D,P)$ (denoted by $(H,F^{\cdot}_{\cdot},G^{\cdot}_{\cdot})$)  with origin ${\bf c}$} 
consists of a vector bundle $H$
on $X$ together with filtrations $F^i$  on the restrictions 
$H_{D_i}$ of $H$ on $D_i$, and furthermore filtrations $G^j$ of the vector spaces $H_{P_j}$. 
The indexing of these filtrations is $F^i_j$ for $c_{D,i} -n \leq  j \leq c_{D,i}$ 
with $F^i_{c_{D,i}}= H|_{D_i}$ and $F^i_{c_{D,i}-n}=0$, and $G^j_k$ for $c_{P,j}-n \leq k \leq c_{P,j}$
with analogous end conditions. 
\end{definition}

Let $H'=\pi^*H$ be the pullback of the vector bundle $H$.
The filtrations $F^j_i$ along the $D'_i$ and $G^j_k$ along the exceptional divisors $E_j$
determine a parabolic structure denoted $\Phi (H,F^{\cdot}_{\cdot},G^{\cdot}_{\cdot})$
over $(X',D'+E)$. By Lemma \ref{surface}, it is automatically locally abelian.

We can use the formula of Theorem \ref{answer} to obtain a formula for the Chern character of 
$\Phi (H,F^{\cdot}_{\cdot},G^{\cdot}_{\cdot})$ 


Consider the push--forward map
$$
\pi_*:CH_.(X')\otimes \Q \lrar CH_.(X)\otimes \Q
$$

We define the Chern character of the exceptionally constant parabolic structure on $X$,
$(H,F^{\cdot}_{\cdot},G^{\cdot}_{\cdot})$, to be 
$$
\m{ch}(H,F^{\cdot}_{\cdot},G^{\cdot}_{\cdot}):= \pi_*\m{ch } \Phi (H,F^{\cdot}_{\cdot},G^{\cdot}_{\cdot}).
$$


\subsection{Parabolic bundles associated to unipotent monodromy at infinity}

Recall that one can associate a parabolic bundle to a logarithmic connection with rational residues,
in a canonical way, such that the weights correspond to the eigenvalues of the residues 
(see \cite{Iy-Si} or \S \ref{reznikov} below). 
In this section, we  point out that one can do something substntially 
different, in the case of nilpotent residues. Suppose $(E,\nabla )$ is a logarithmic connection on $X$, with singularities along a normal-crossings
divisor $D = D_1+ \ldots + D_m$, such that the residue $\eta_i$ of $\nabla$ are nilpotent, for each $i=1,...,m$. In other words, $(E,\nabla )$ is the Deligne
extension of a flat bundle with unipotent monodromy at infinity. 

In this case, we still have some different natural filtrations along divisor components, but the eigenvalues of 
the residue are zero so there is no canonical choice of weights. Instead, define some
characteristic numbers by arbitrarily assigning weights to these filtrations. Assume that $X$ is
a surface here, so that the resulting parabolic structures will automatically be locally abelian.
It seems to be an interesting question to determine when the locally abelian condition holds for
these kinds of filtrations in the case of dimension $\geq 3$. 

Consider the {\em Image filtration} on the restriction $E_{D_i}$ of $E$ to a divisor component:
$$
E_{D_i}=F^i_0\supset F^i_1\supset ...\supset F^i_{l_i-1}\supset F^i_{l_i+1}=0
$$
where
$$
F^i_j:= \m{image }(\eta_i^j: E_{D_i}\lrar E_{D_i}),
$$
$\eta_i^j:= \eta_i\circ \eta_i\circ...\circ \eta_i$ ($j$-times) and $l_i+1$ is the order of $\eta_i$.
 
Alternatively, we can consider the {\em Kernel filtration} induced by the kernels of the operator $\eta_i$:
write
$$
F^i_j:= \m{kernel }(\eta_i^{l_i+1-j}: E_{D_i}\lrar E_{D_i}).
$$ 
   
Mixing these two filtrations gives rise to the {\em monodromy weight filtration } $\{W_l\}$ defined by Deligne \cite{Deligne3}.
This is an increasing filtration
$$
\{0\}\subset W_0\subset W_1\subset...\subset W_{2l_i}=E_{D_i}
$$
uniquely determined by the conditions:

$\bullet\,\,\eta_i(W_l)\subset W_{l-2}
$

$
\bullet\,\,\m{ the induced map }\eta_i^l: \m{Gr}_{k+l}(W_*)\rar \m{Gr}_{k-l}(W_*) 
$
$ \m{ is an isomorphism for each }l.
$

Here $\m{Gr}_l(W_*):= W_l/W_{l-1}$.

Explicitly, the filtration is defined by induction as follows:
let 
$$
W_0=\m{image}(\eta_i^{l_l}) \m{ and } W_{2l_i-1}= \m{ ker}(\eta_i^{l_i}).
$$
Now fix some $l< l_i+1$; if
$$
0\subset W_{l-1}\subset W_{2l_i-l}\subset W_{2l_i}=E_{D_i}
$$
has already been defined in such a way that 
$$
\eta_i^{l_i-l+1}(W_{2l_i-l})\subset W_{l-1}
$$
then we define
$$
W_l/W_{l-1}=\m{image} (\eta_i^{l_i-l}:{W_{2l_i-l}/W_{l-1}} \lrar {W_{2l_i-l}/W_{l-1}}    )
$$
and $W_l,W_{2l_i-l-1}$ to be the corresponding inverse images. Notice that 
$$
W_l/W_{l-1}\subset W_{2l_i-l-1}/W_{l-1}
$$ 
so that $W_l\subset W_{2l_i-l-1}$. Clearly, 
$\eta_i^{l_i-1}(W_{2l_i-l-1})\subset W_l$, so that the induction hypothesis is satisfied.

\begin{lemma}\label{parsurface}
Suppose $X$ is a surface. 
Consider the {\em Image } or the {\em Kernel} or the {\em monodromy weight filtrations} considered above, on the 
restrictions $E_{D_i}$ of $E$ to the divisor components. We can associate a 
locally abelian parabolic bundle  on $(X,D)$ with respect to $(E_U,\nabla_U)$ 
together with either of these filtrations by assigning aribitrary weights.
\end{lemma}
\begin{proof}
By Lemma \ref{surface}, the parabolic structure defined by the filtrations is automatically locally abelian. 
\end{proof}


\subsection{Examples arising from families}

Suppose $\pi:X\lrar S$ is a semi-stable family of projective varieties such that
$\pi_U:X_U\lrar U$ is a smooth morphism, for some open subvariety $U\subset S$ and $D:=S-U$ is a normal crossing divisor. Let $d$ be the relative dimension of $X\lrar S$.

In this situation, the Gauss--Manin bundles $\cH^l:=R^l\pi_*(\Omega_{X/S}^\bullet(\pi^{-1}D))$
for $0\leq l\leq 2d$, are equipped with a logarithmic flat connection. Furthermore, the local monodromies are unipotent and $\cH^l$ is the Deligne extension of the restriction $\cH^l_U$ (see \cite{Steenbrink}).
Let $\eta_i$ be the residue transformations along the divisor components $D_i$. Unipotency of 
the monodromy operators implies nilpotency of $\eta_i$ and the order of nilpotency is at most $l+1$ (see \cite{Landman}).
In particular, the length of the Image and the Kernel filtrations in the previous 
subsection is at most $l+1$ and the monodromy weight filtration is of length at most $2l+1$. 
We make an explicit computation of the Chern character of the associated locally abelian parabolic bundle in the following case:

Suppose $S$ is a surface and $X\lrar S$ is a semi-stable family of abelian varieties.
We consider the Gauss-Manin system $\cH^1$ of weight one on $S$. 
For simplicity assume that $D$ is a smooth irreducible divisor.
Then the residue transformation $\eta$ has order of nilpotency  two and in this case the \textit{monodromy weight filtration} is written as
$$
{\cH^1}_{|D}=W_2\supset W_1 \supset W_0 \supset W_{-1}=0.
$$
Here $W_1= \m{kernel} (\eta)$ and $W_0= \m{image}(\eta)$. The graded pieces
$$
\textbf{gr}_m:= \f{W_m}{W_{m-1}} 
$$
carries a polarized pure Hodge structure of weight $m$ (see \cite{Schmid}).
Also, the graded piece of weight two is isomorphic to the piece of weight zero, by the monodromy operator $N$
(in \cite{Schmid}, $N$ polarizes the mixed Hodge structures). 

By Lemma \ref{parsurface}, we can associate a locally abelian parabolic bundle $F$ on $S$ corresponding to $\{W_.\}$, with arbitrary weights  
$(\alpha _0, \alpha _1, \alpha _2)$
with $-1 < \alpha _0 < \alpha _1 < \alpha _2 \leq 0$.

\begin{lemma}\label{weightone}
Suppose $X\lrar S$ is a semi-stable family of abelian varieties of genus $g$. Let $g_i$ denote the rank of $\textbf{gr}_i$ for $i=0,1,2$,
thus $g= g_0+g_1+g_2$ and $g_0=g_2$. With notations as above, assigning weights $(\alpha _0, \alpha _1, \alpha _2)$,
the Chern character of the locally abelian parabolic bundle $F$ is given by the formula
$$
{\rm ch}(F) = \sum _{i=0}^2 g_ie ^{-\alpha _iD} \in CH^{\cdot}(S)_{\Q}.
$$
In other words it is Chow-equivalent to a direct sum of parabolic line bundles. 
\end{lemma}
\begin{proof}
Let $k:D\hookrightarrow X$ denote the inclusion. Suppose $A$ is a rank $r$ bundle along $D$ whose Chern character is $r\in CH^0(D)_{\Q}$.
Then, the sheaf $k_{\ast}(A)$ on $X$ has Chern character given by a Riemann-Roch formula. This formula depends only on the Chern character
of $A$ on $D$, in particular it is $r$ times the value for the case $A=\cO _D$. In that case we can use the exact sequence
$$
0\rightarrow \cO (-D) \rightarrow \cO \rightarrow \cO _D \rightarrow 0
$$
to conclude that the Chern character of $k_{\ast}(A)$ is $r(1-e ^{-D})$. 

Turn now to the situation of the lemma. By \cite{vanderGeer} or \cite{EsnaultViehweg}, we have 
$$
\m{ch}(E) = g \in CH^0(S)_{\Q}
$$
and similarly for $\m{ch}(\textbf{gr}_1)$ which corresponds to a family of abelian varieties along $D$, we get
$$
\m{ch}(\textbf{gr}_1) =g_1\in CH^0(D)_{\Q} .
$$
Clearly, $\m{ch}(\textbf{gr}_0)= g_0 \in CH^0(D)_{\Q}$, thus 
$\m{ch}(\textbf{gr}_2)= g_2 \in CH^0(D)_{\Q}$ by the isomorphism between the weight two and weight zero piece given by the monodromy operator. 
Plugging these into the formula \eqref{onecomponent2}  of \S \ref{subsec.-singlesmooth} and using the previous paragraph for
the Chern characters of $k_{\ast}(A)$ we get the formula 
$$
{\rm ch}(F) = e^Dg - \sum _{i=0}^2 \frac{e^D - e ^{-\alpha _iD}}{(1-e^{-D} )}g_i(1-e ^{-D}) \in CH^{\cdot}(S)_{\Q}.
$$ 
Simplifying with $g=g_0+g_1+g_2$ gives the stated formula. 
\end{proof}


\section{Extended Reznikov theory for finite order monodromy at infinity}
\label{reznikov}

Suppose $U$ is a nonsingular variety defined over the complex numbers. 
Consider a nonsingular compactification $X$ of $U$ such that $D:=X-U$ is a normal crossing divisor. 
Suppose $(E_U,\nabla_U)$ is a bundle with a flat connection on $U$. 
Consider the canonical extension $(E,\nabla)$ of $(E_U,\nabla_U)$ on $X$ (see \cite{Deligne}). 
Here $\nabla$ is a logarithmic connection on $E$, i.e.,
$$
\nabla:E\lrar E \otimes \Omega_X(\m{log} D)
$$
is a $\comx$-linear map and satisfies the Leibnitz rule.
Flatness implies that $\nabla\circ \nabla =0$.

Consider the sequence induced by the Poincar\'e residue map

$$
E\lrar E\otimes \Omega_X(\m{log}D) \sta{res}{\lrar} E\otimes \cO_D.
$$
This induces an operator

$$
\eta_i: E_{D_i}\lrar (E\otimes \Omega_X(\m{log}D))_{|D_i} \sta{res}{\lrar} E_{D_i}
$$
called the residue transformation along the divisor component $D_i$ and $\eta_i\,\in\, \m{End}(E_{D_i})$.   

\begin{definition}
\label{rationalresidues}
We say that $(E,\nabla )$ {\em has rational residues} if
 the eigenvalues of the residual transformations $\eta _i$ above are rational numbers.
\end{definition}

This is equivalent to saying that the local monodromy transformations around the divisor components $D_i$ of $D$ are quasi-unipotent.

If $\al_i$ are the rational residues then \cite{Deligne}
$$
e^{2\pi i\al_i}=\m{eigenvalues of the local monodromy}.
$$

Suppose the residues of $(E,\nabla)$ are non-zero and rational.
In \cite[Lemma 3.3]{Iy-Si}, a locally abelian parabolic bundle $\cE$ on $(X,D)$ was associated to $(E,\nabla)$. 
In fact, $\cE$ was associated to the flat connection $(E_U,\nabla_U)$ on $U$ and the constituent bundles were defined, 
using a construction due to Deligne-Manin. If we choose the extension $(E,\nabla)$ on $X$ such that the 
rational residues lie in the interval $[0,1)$ then the weights are precisely the inverse of the rational 
residues. In other words, if $0\leq -\al_i^1<-\al_i^2<...<-\al_i^{n_i}< 1$ are the rational residues along $D_i$ 
then the weights are $\al_i^1>\al_i^2>...>\al_i^{n_i}$ along $D_i$.

\subsection{Residues are rational and semisimple}
Suppose that the residues are rational and furthermore on the associated-graded of the parabolic structure, the
residue of the connection induces a semisimple operator. In this case, the monodromy transformations
of the corresponding local system are semisimple with eigenvalues which are roots of unity, thus they are
of finite order. If $n$ denotes a common denominator for the rational residues of the connection
(and hence for the corresponding parabolic weights) then the monodromy transformations have order $n$.
This implies that the connection extends to a flat connection on
the DM-stack $Z:= X \langle \frac{D_1}{n},\ldots , \frac{D_m}{n}\rangle$. Conversely any flat connection on the DM-stack
$Z$ gives rise to a connection on $U$ with semisimple and rational residues.

The locally abelian parabolic bundle on $X$ corresponds to the vector bundle on $Z$ underlying the
flat bundle as extended over $Z$. Indeed, when the monodromy transformations have order $n$, the monodromy around the divisor at
infinity in $Z$ is trivial, and in this case the Deligne canonical extension is the vector bundle underlying the extended flat bundle.
By \cite{Iy-Si} the Deligne canonical extension over $Z$ is the vector bundle corresponding to the parabolic bundle on $X$.  


\subsection{Reznikov's theory in the case of rational semisimple residues}

The theory of secondary characteristic classes works equally well on the DM-stack $Z$. In particular, we can define
the rational Deligne cohomology 
$$
H^{2p}_{\cD}(Z, \Q (p)):= \hH ^{2p}(Z^{\rm an}; \Q(p) \rightarrow \Omega ^0_{Z^{\rm an}} \rightarrow \ldots 
\rightarrow \Omega ^p_{Z^{\rm an}}),
$$
and also the cohomology 
$$
H^{2p-1}(Z, \comx / \Q ) = \hH ^{2p}(Z^{\rm an}; \Q \rightarrow \Omega ^{\cdot}_{Z^{\rm an}}).
$$
Dividing by the Hodge filtration provides a map 
\begin{equation}
\label{dividehodge}
H^{2p-1}(Z, \comx / \Q ) \rightarrow H^{2p}_{\cD}(Z, \Q (p)).
\end{equation}
On the other hand, the Deligne cycle class map from Chow groups to Deligne cohomology is a map
\begin{equation}
\label{delignecycle}
CH^p(Z)_{\Q} \rightarrow H^{2p}_{\cD}(Z, \Q (p)).
\end{equation}
If $E$ is a vector bundle on $Z$ then its Chern character in $CH^{\cdot}(Z)_{\Q}$ maps to its Chern character in 
$\oplus _pH^{2p}_{\cD}(Z, \Q (p))$. 

\begin{lemma}
\label{deligneisoXZ}
Pullback for the map $Z\rightarrow X$ gives an isomorphism of 
Deligne cohomology groups
$$
H^{2p}_{\cD}(X, \Q (p)).
\stackrel{\cong}{\longrightarrow} H^{2p}_{\cD}(Z, \Q (p))
$$
compatible with the isomorphism of rational Chow groups and the map \eqref{delignecycle}. 
It also induces an isomorphism 
$$
H^{2p-1}(X, \comx / \Q ) \stackrel{\cong}{\longrightarrow}H^{2p-1}(Z, \comx / \Q ) 
$$
and this is compatible with the projection \eqref{dividehodge}.
\end{lemma}

Suppose $F$ is a locally abelian parabolic bundle on $X$. Define the Chern character of $F$ in Deligne cohomology
of $X$ by using the formula of Theorem \ref{answer} and taking the 
Chern characters of the pieces $F_{[a_1,\ldots ,a_m]}$ in
the Deligne cohomology of $X$. Thus 
\begin{equation}
\label{delignechdef}
{\rm ch}^{\cD}(F) := \frac{\sum _{a_1=0}^{n-1}\cdots \sum _{a_m=0}^{n-1} e^{-\sum _{i=1}^m \f{a_i}{n}c_1^{\cD}(D_i)}
{\rm ch}^{\cD}(F_{[a_1,\ldots , a_m]})}
{\sum _{a_1=0}^{n-1}\cdots \sum _{a_m=0}^{n-1} e^{-\sum _{i=1}^m \f{a_i}{n}c_1^{\cD}(D_i)} }.
\end{equation}
The products are taken with the product structure of Deligne cohomology which is compatible with the
intersection product in Chow groups \cite{Es-Vi2}. 

\begin{corollary}
\label{answerdelcoh}
Suppose $F$ is a locally abelian parabolic bundle on $X$ with $n$ as common denominator for the
rational weights, corresponding to a vector bundle $E$ on $Z$. 
Then ${\rm ch}^{\cD}(F)$ as given by the above formula \eqref{delignechdef}, pulls back to ${\rm ch}^{\cD}(E)$ on $Z$ 
via the isomorphism of Lemma \ref{deligneisoXZ}.
\end{corollary}
\begin{proof}
By Theorem \ref{answer} this is the case for the Chern character in Chow groups, and we have the
compatibility of the isomorphism of Lemma \ref{deligneisoXZ} with the projection  \eqref{delignecycle}. 
\end{proof}

Now, go back to the situation where $(E_U, \nabla _U)$ is a flat bundle on $U$ with rational and semisimple
residues. It extends to a flat bundle 
$(E,\nabla )$ on $Z$ and also the local system $L_U$ on $U$ extends 
to a local system $L$ on $Z$. 

Consider a Kawamata cover (see \cite{Kaw})
$$
f:Y\lrar X
$$
so that $Y$ is a smooth projective variety. Then there is a factorization
$$
Y\sta{h}{\lrar} Z \sta{\pi}{\lrar} X
$$
such that $f=\pi\circ h$. The flat connection on $Z$ pulls back to a flat connection $(E_Y,\nabla_Y)$ on $Y$.
Thus, Esnault's theory of secondary classes for flat bundles \cite{Esnault} gives a class 
$\widehat{c}_p (L)\in H^{2p-1}(Z, \comx / \Q )$. By \cite{Esnault}, this class projects 
under the map \eqref{dividehodge} to the Deligne Chern class $c_p^{\cD}(E)$ for the vector bundle $E$ on $Z$.

\begin{proposition}\label{parReznikov}
Reznikov's result on the vanishing of the rational secondary classes works equally well over a smooth projective
DM-stack. Thus, with the above notations  $\widehat{c}_p (L)=0$ in $H^{2p-1}(Z, \comx / \Q )$, for $p\geq 2$. 
\end{proposition}
\begin{proof}
Either of Reznikov's proofs of \cite{Reznikov} work equally well over the DM-stack $Z$. 
Alternatively, we can reduce to the utilisation of \cite{Reznikov} on the finite cover $Y$ as follows:
by Reznikov's theorem the secondary classes of $(E_Y,\nabla_Y)$ are 
trivial in the $\comx/\Q$-cohomology in degrees $\geq 3$ of $Y$. 
The map $Y\rightarrow Z$ induces an injection $H^i(Z,V)\hookrightarrow H^i(Y,V)$ for any $\Q$-vector space
$V$, in particular $V=\comx / \Q$. 
This implies that the secondary classes  vanish on $Z$
$$
\widehat{c}_p (L)=0\, \in\, H^{2p-1}(Z, \comx / \Q )
$$
for $p\geq 2$.

\end{proof}

Combining with our formula of Theorem \ref{answer} we obtain a formula for an element of the Deligne cohomology
over the compactification $X$ of $U$ which vanishes by Reznikov's theorem. 

\begin{corollary}
\label{parabolicreznikov}
Suppose $(E_U,\nabla _U)$ is a flat bundle on $U$ with rational and semisimple residues, or equivalently
the monodromy transformations at infinity are of finite order. Let $F$ denote the corresponding locally 
abelian parabolic bundle. Define the Deligne Chern character ${\rm ch}^{\cD}(F)$ on $X$ by the formula
\eqref{delignechdef}. Then the Deligne Chern classes ${\rm c}^{\cD}_p(F)$ in all degrees $\geq 2$ vanish.
\end{corollary}
\begin{proof}
This follows from Corollary \ref{answerdelcoh} and Proposition \ref{parReznikov}.
\end{proof}



\end{document}